\documentclass[epic,11pt]{article}
\usepackage{amsfonts,amsmath,amscd,amssymb,amsthm}
\usepackage{fullpage}
\usepackage[alphabetic,initials]{amsrefs}
\usepackage[all]{xy}
\usepackage{mathrsfs}
\input xy%
\usepackage{fullpage}
\def\convf{\hbox{\space \raise-2mm\hbox{$\textstyle      \bigotimes \atop \scriptstyle \omega$} \space}}
\def\0{{\bar 0}}
\def\1{{\bar 1}}
\def\mo{{\operatorname{Mor}}}

\def\Z{{\mathbb Z}}

\def\B{{\mathcal B}}

\def\Vc{{V^c}}

\def\N{{\mathbb N}}

\def\Max{{\operatorname{Max}\;}}
\def\id {{\operatorname{id}}}
\def\reg {{\operatorname{reg}}}

\def\Spec  {{\operatorname{Spec  }\;}}
\def\ev  {{\operatorname{ev  }}}
\def\ad  {{\operatorname{ad  }\;}}
\def\Ad  {{\operatorname{Ad  }}}
\def\Rad{\operatorname{Rad\;}}
\def\rad{\operatorname{rad\;}}
\def\span{\operatorname{span\;}}
\def\Ker {{\operatorname{Ker}\;}}
\def\Lie {{\operatorname{Lie}\;}}

\def\heart {{\operatorname{heart}}}

\def\Id{{\operatorname{Id}}} 

\def\atyp{{\operatorname{atyp}}}
\def\iso{{\operatorname{iso}}}

\def\Hom{{\operatorname{Hom}}}

\newcommand{\ttE}{\mathtt{E}}
\newcommand{\ttF}{\mathtt{F}}

\newcommand{\ttk}{\mathtt{k}}
\newcommand{\tte}{\mathtt{e}}

\newcommand{\ttT}{\mathtt{T}}

\newcommand{\itemi}{\item[{{\rm(i)}}]}
\newcommand{\itemii}{\item[{{\rm(ii)}}]}
\newcommand{\itemiii}{\item[{{\rm(iii)}}]}
\newcommand{\itemiv}{\item[{{\rm(iv)}}]}
\newcommand{\itemv}{\item[{{\rm(v)}}]}

\newcommand{\itema}{\item[{{\rm$($a$)$}}]}
\newcommand{\itemb}{\item[{{\rm$($b$)$}}]}
\newcommand{\itemc}{\item[{{\rm$($c$)$}}]}
\newcommand{\itemd}{\item[{{\rm$($d$)$}}]}
\newcommand{\iteme}{\item[{{\rm$($e$)$}}]}

\newcommand{\itemo}{\item[{}]}
\newcommand{\noi}{\noindent}

\newcommand{\ga}{\alpha}
\newcommand{\gb}{\beta}
\newcommand{\gc}{\gamma}

\newcommand{\Gc}{\Gamma}
\newcommand{\Gl}{\Lambda}
\newcommand{\Gd}{\Delta}

\newcommand{\gd}{\delta}

\newcommand{\gs}{\sigma}
\newcommand{\gS}{\Sigma}

\newcommand{\gO}{\Omega}

\newcommand{\gt}{\tau}

\newcommand{\gl}{\lambda}

\newcommand{\gr}{\rho}

\newcommand{\gep}{\epsilon}

\newcommand{\op}{\oplus}
\newcommand{\A}{\mathcal A}

\def\res{{\operatorname{res}}}

\def\Im{{\operatorname{Im}\;}}

\newcommand{\ot}{\otimes}

\newcommand{\fg}{\mathfrak{g}}\newcommand{\fgl}{\mathfrak{gl}}
\newcommand{\fsl}{\mathfrak{sl}}\newcommand{\fpsl}{\mathfrak{psl}}\newcommand{\osp}{\mathfrak{osp}}

\newcommand{\fh}{\mathfrak{h}}

\newcommand{\fb}{\mathfrak{b}}

\newcommand{\fa}{\mathfrak{a}}

\newcommand{\fp}{\mathfrak{p}}
\newcommand{\fq}{\mathfrak{q}}
\newcommand{\pq}{\mathfrak{pq}}
\newcommand{\psq}{\mathfrak{psq}}

\newcommand{\fA}{\mathfrak{A}}
\newcommand{\fB}{\mathfrak{B}}
\newcommand{\fG}{\mathfrak{G}}

\newcommand{\sfp}{{\mathfrak{sp}}}
\newcommand{\sfq}{\small{\mathfrak{sq}}}

\newcommand{\ff}{\footnote}
\newfont{\eufm}{eufm10 scaled\magstep1}

\newcommand{\ci}{\circ} \newcommand{\ti}{\times}

\newcommand{\bca}{\bigcap}
\newcommand{\bcu}{\bigcup}
\newcommand{\bop}{\bigoplus}
\newcommand{\bsk}{\backslash}

\newcommand{\cO}{\mathcal{O}}

\newcommand{\cI}{\mathcal{I}}

\newcommand{\cA}{\mathcal{A}}

\newcommand{\cF}{\mathcal{F}}

\newcommand{\cG}{\mathcal{G}}

\newcommand{\cS}{\mathcal{S}}
\newcommand{\cT}{\mathbb{T}}

\newcommand{\cV}{\mathcal{V}}

\newcommand{\cW}{\mathfrak{W}}

\newcommand{\ey}{\end{eqnarray}}
\newcommand{\by}{\begin{eqnarray}}
\newcommand{\nn}{\nonumber}

\newcommand{\bco}{\begin{conjecture}}
\newcommand{\ba}{\begin{alg}}
\newcommand{\ea}{\end{alg}}
\newcommand{\eco}{\end{conjecture}}
\newcommand{\bpf}{\begin{proof}}
\newcommand{\epf}{\end{proof}}
\newcommand{\bt}{\begin{theorem}}

\newcommand{\et}{\end{theorem}}

\newcommand{\br}{\begin{rem}}
\newcommand{\er}{\end{rem}}
\newcommand{\brs}{\begin{rems}}
\newcommand{\ers}{\end{rems}}
\newcommand{\bi}{\begin{itemize}}
\newcommand{\ei}{\end{itemize}}
\newcommand{\bl}{\begin{lemma}}
\newcommand{\bsul}{\begin{sublemma}}
\newcommand{\esul}{\end{sublemma}}
\newcommand{\bp}{\begin{proposition}}
\newcommand{\be}{\begin{equation}}
\newcommand{\bc}{\begin{corollary}}
\newcommand{\bexs}{\begin{examples}}
\newcommand{\eexs}{\end{examples}}
\newcommand{\bexa}{\begin{example}}
\newcommand{\eexa}{\end{example}}
\newcommand{\bex}{\begin{exercise}}
\newcommand{\eex}{\end{exercise}}
\newcommand{\btab}{\begin{tab}}
\newcommand{\etab}{\end{tab}}
\newcommand{\el}{\end{lemma}}
\newcommand{\ep}{\end{proposition}}
\newcommand{\ee}{\end{equation}}
\newcommand{\ec}{\end{corollary}}
\newcommand{\Bc}{\begin{center}}
\newcommand{\Ec}{\end{center}}
\newcommand{\bh}{\begin{hyp}}
\newcommand{\eh}{\end{hyp}}
\newcommand{\bhs}{\begin{hyps}}
\newcommand{\ehs}{\end{hyps}}
\newcommand{\bd}{\begin{dfn}}
\newcommand{\ed}{\end{dfn}}

\newcommand{\bn}{\begin{notn}}
\newcommand{\en}{\end{notn}}

\begin{document}
\title{Table of Contents}

\newtheorem*{bend}{Dangerous Bend}

\newtheorem{thm}{Theorem}[section]
\newtheorem{hyp}[thm]{Hypothesis}
 \newtheorem{hyps}[thm]{Hypotheses}
\newtheorem{notn}[thm]{Notation}

  \newtheorem{rems}[thm]{Remarks}

\newtheorem{conjecture}[thm]{Conjecture}
\newtheorem{theorem}[thm]{Theorem}
\newtheorem{theorem a}[thm]{Theorem A}
\newtheorem{example}[thm]{Example}
\newtheorem{examples}[thm]{Examples}
\newtheorem{corollary}[thm]{Corollary}
\newtheorem{rem}[thm]{Remark}
\newtheorem{lemma}[thm]{Lemma}
\newtheorem{sublemma}[thm]{Sublemma}
\newtheorem{cor}[thm]{Corollary}
\newtheorem{proposition}[thm]{Proposition}
\newtheorem{exs}[thm]{Examples}
\newtheorem{ex}[thm]{Example}
\newtheorem{exercise}[thm]{Exercise}
\numberwithin{equation}{section}%
\setcounter{part}{0}
\newcommand{\drar}{\rightarrow}
\newcommand{\lra}{\longrightarrow}
\newcommand{\rra}{\longleftarrow}
\newcommand{\dra}{\Rightarrow}
\newcommand{\dla}{\Leftarrow}
\newcommand{\rl}{\longleftrightarrow}

\newtheorem{Thm}{Main Theorem}


\newtheorem*{thm*}{Theorem}
\newtheorem{lem}[thm]{Lemma}
\newtheorem*{lem*}{Lemma}
\newtheorem*{prop*}{Proposition}
\newtheorem*{cor*}{Corollary}
\newtheorem{dfn}[thm]{Definition}
\newtheorem*{defn*}{Definition}
\newtheorem{notadefn}[thm]{Notation and Definition}
\newtheorem*{notadefn*}{Notation and Definition}
\newtheorem{nota}[thm]{Notation}
\newtheorem*{nota*}{Notation}
\newtheorem{note}[thm]{Remark}
\newtheorem*{note*}{Remark}
\newtheorem*{notes*}{Remarks}
\newtheorem{hypo}[thm]{Hypothesis}
\newtheorem*{ex*}{Example}
\newtheorem{prob}[thm]{Problems}
\newtheorem{conj}[thm]{Conjecture}

\title{On the geometry of some algebras related to the Weyl groupoid.}
\author{Ian M. Musson
 \\Department of Mathematical Sciences\\
University of Wisconsin-Milwaukee\\ email: {\tt
musson@uwm.edu}}
\maketitle
\begin{abstract} Let $\ttk$ be an algebraically closed field of characteristic zero.  Let $\fg$ be a finite dimensional classical simple Lie superalgebra over $\ttk$ or $\fgl(m,n)$. In the case that $\fg$ is a Kac-Moody algebra of finite  type with set of roots $\Gd$, Sergeev and Veselov introduced the 
Weyl groupoid $\mathfrak{W}=\mathfrak{W}(\Gd)$, which has significant connections with the representation theory of $\fg$. Let $\fh$, $W$ and $Z(\fg)$  be a Cartan subalgebra of 
$\fg_0$, the Weyl group of $\fg_0$ and  the center of $U(\fg)$ respectively. Also let $G$ be a Lie  supergroup  with Lie $G =\fg$. There are several important
commutative algebras related to $\mathfrak{W}$.  Namely 
\bi \item The image $I(\fh)$ of the injective Harish-Chandra map $Z(\fg)\lra S(\fh)^W$.
\item The supercharacter $\Z$-algebras $J(\fg)$ and $J(G)$ of finite dimensional representations of $\fg$ and $G$.
\ei
Let $\A = \A(\fg)$ be denote either $I(\fh)$ or $J(G) \ot_{\Z}\ttk$. 
The purpose of this  paper 
 is to investigate the algebraic geometry of $\A.$  In many cases, the algebra $\A$ satisfies    
the Nullstellensatz. This gives a bijection between radical ideals in $\A$ and superalgebraic sets (zero loci of such ideals).  Any  superalgebraic set
is uniquely a finite union of irreducible superalgebraic  components. In the  
Kac-Moody case, we describe the smallest superalgebraic set containing a given (Zariski) closed set, and show that the superalgebraic sets are exactly the closed sets that are unions of groupoid orbits.
\end{abstract}
\noi 
\ref{se1}. {Introduction.}\\ \ref{se2}. {General Results.}\\ \ref{se3}. {Applications.}\\  \ref{sn}. {The Nullstellensatz.}\\ \ref{sas}. {Weyl Groupoids and their Actions.}\\\ref{csc}. {The Geometry of Superalgebraic Sets.}
\section{Introduction} \label{se1}
\subsection{Background. } \label{Br}
\subsubsection{The classical Lie superalgebras. } \label{Cls}
We are interested in the following {\it classical} Lie superalgebras $\fg$.  All have $\fg_0$ reductive and $\fg_0\neq \fg$. 
\bi \itema The KM \ff{Kac-Moody Lie superalgebras of finite type. These algebras are also known as contragredient Lie superalgebras, \cite{Kac1}, \cite{M} Chapter 5.} algebras 
$\fgl(m|n),$ 
$    \osp(r|2n), $ 
$\fsl(m|n)$ with $m\neq n,$  
$    D(2|1;\ga),$ 
$    F(4),$ $G(3),$ 
\itemb The strange Lie superalgebras $\fp(n)$, $\fq(n)$ and their simple relatives $\sfp(n)$, $\psq(n)$.
\ei
There are some conditions on $r, m,n\in\N$ and $\ga\in\ttk$ which we omit.  
 The last three algebras listed in (a) are called {\it exceptional}. 
The algebras listed in (a) with $ \fgl(m|n)$ replaced  by  $\fpsl(n|n)$, together with the simple algebras from (b) are  precisely the  finite dimensional simple Lie superalgebras $\fg$ with  $\fg_0$ reductive and $\fg_0\neq \fg$. 
\noi 

\subsubsection{Algebras related to the Weyl groupoid. } \label{Awg}
 The Weyl groupoid $\mathfrak{W}(\Gd)$ of Sergeev and Veselov, \cite{SV2} 
  has 
become a fundamental tool in the  representation theory  of a KM algebra $\fg$. Here $\Gd$ is the root system of  $\fg$. This groupoid was introduced in relation to the algebras 
$J(\fg)$ and $J(G)$ listed in the abstract.  The
Lie supergroups $G$ we  consider are defined  in Lemma \ref{abc}.  In particular the supergroups $P(n)$ and $Q(n)$ have
 $\Lie P(n)=\fp(n)$ and $\Lie Q(n)=\fq(n)$. 
Weyl groupoids for   $ P(n)$ and $Q(n) $ were defined in  \cite{IRS} Section 5.4
and  \cite{Re}, Section 5 
 respectively.  
\\ \\
For  a KM algebra $\fg $, an important role is played by a non-degenerate, symmetric  bilinear form   $(\;,\;)$ on  $\fg $.  This can be constructed in a uniform way,  \cite{M} Theorem 5.4.1, Remark 5.4.2.  
The restriction of this form to $\fh$ remains non-degenerate 
and this yields an isomorphism $\fh \lra \fh^*$, denoted  $h_\ga\mapsto \ga$. Thus we obtain a non-degenerate, symmetric  bilinear form on
$\fh^*$  which we also denote by $(\;,\;)$.  This latter form 
is invariant 
under the Weyl group $W$ of $\fg_0$, 
and we  can normalize the whole construction by requiring that its Gram matrix is the (symmetrized) Cartan matrix $(\ga_i,\ga_j)$, where 
 $\{\ga_i\}$ is  the distinguished basis of simple roots for  $\fg$, 
\cite{M} Section 3.2. 
The issue with $ \fpsl(n|n)$ is
 that  any set of simple roots 
$\{ \alpha_i  \}$ is linearly dependent, and the  matrix $A= (  \alpha_i,  \alpha_{j})$ is  singular.  Thus there is no bilinear form as above. The usual solution to this problem is to take a minimal realization of $A$, leading  to the construction of $\fgl(n|n),$
which  has $\fpsl(n|n)$ as a subquotient. 
For  $\alpha \in 
\Gd_{iso}$, the set of isotropic roots of $\fg$,  set \be \label{rtal}  \Pi_{\ga}=
 \{ \lambda \in   \fh^* | (\lambda,\ga) =0\}.
\ee   
The other algebra of interest $I(\fh)$ is isomorphic to   
$Z(\fg) $ via the  Harish-Chandra map.
This  algebra may be defined as 
\be \label{ltal}  I(\fh)= \{f\in S(\mathfrak{h})^W| f(\lambda) = f(\lambda + t \alpha) \mbox{ for all }\alpha \in 
\Gd_{iso}, \gl \in \Pi_{\ga} \mbox{ and } t \in \ttk\}.\ee Several other definitions of $I(\fh)$ can be given based on \cite{M}, Lemma 12.1.1.
By work of Gorelik and Kac \cite{Gk}, \cite{Kac4}, the  Harish-Chandra map yields an isomorphism $Z(\fg)\lra I(\fh)$, see  also \cite{M}, Theorem 13.1.1.  
For a description of $Z(\fg)$ when  $\fg=\fp(n)$ and $\fq(n)$, see \cite{Gor},  \cite{Gq} respectively.
\\ \\
There is a description of $J(G)$ parallel to \eqref{ltal}. 
Let $\mathbb{T}$ be a maximal torus  in $G_0$, 
and 
$X(\mathbb{T}) = \Hom(\mathbb{T},\ttk^*) $ (resp. 
$Y(\mathbb{T}) = \Hom(\ttk^*, \mathbb{T})$)  the group of 
characters (resp. one-parameter subgroups) of  $\mathbb{T}$. 
Let  
 $D_\ga$ be the derivation of  $\Z[X(\cT)] $ given by  
\be \label{pvp} D_\ga(\tte^\gb) = (\ga, \gb) \tte^\gb.\nn\ee
Then in the KM case, we have by \cite{SV2} Equation (1), 
\be \label{yta}  J(G)= \{f \in \Z[X(\cT)]^W| D_\ga f \in (\tte^\ga -1)  \mbox{ for all }\alpha \in 
\Gd_{iso}\}.\ee 
Various conditions equivalent to that defining $J(G)$ are given in 
 \cite{M22} Section 3.
\subsubsection{Continuous Weyl groupoids. } \label{Cwg}
In the theory of algebraic groups acting on varieties, a key role is played by closed orbits.  However the orbits of 
  $\mathfrak{W}$ on  $\fh^*$ or  $ \mathbb{T}$  are not closed unless they are finite, see Remark \ref{sgs1}, and for  this reason in Section \ref{sas}  
we introduce the continuous   Weyl groupoids $\cW^c$  and $\cW_*^c$. 
The action of a groupoid $\fG$ on an affine variety $X$ and  the invariant ring $\cO(X)^\fG$ are defined 
in Section \ref{sas}. 
  In Proposition \ref{qag}, we show there are actions of  $\cW^c$  and $\cW_*^c$ 
on $\fh^* $  and $ \mathbb{T}$ respectively such that 
$$I(\fh) =S(\fh)^{\cW^c} \mbox{ and }  J(G) \ot_{\Z}\ttk =\cO(\mathbb{T})^{\cW_*^c}.$$
\bn \label{bn}{\rm We assume one of the following holds
\bi \itemi $X=\mathbb{\fh^*}, \; \fG = \mathfrak{W}^c,\;  \B = \cO(X)^W= S(\fh)^W$ and  $\A = I(\fh)$.
\itemii $X=\mathbb{T}, \; \fG = \mathfrak{W}_*^c,\;  \B = \cO(X)^W$ and  $\A = J(G) \ot_{\Z}\ttk$. 
\ei 
}
\en \noi
\bexa{\rm If $\fg=\fgl(m|n)$ the algebras $I(\fh)$ and  $J(G) \ot_{\Z}\ttk$ are the algebras of 
 supersymmetric polynomials and Laurent supersymmetric polynomials over $\ttk$ respectively.  These algebras have very beautiful combinatorics, see \cite{M} Chapter12, \cite{PT}  and  \cite{Se}. Suprisingly though, the algebras of (Laurent) supersymmetric polynomials have received little attention from a geometric point of view, perhaps because they are not Noetherian.\\ \\ In the simplest case where $m=n=1$, $\cW$ is the smallest groupoid that is not a disjoint union of groups.  Thus $\cW$  has two objects, and two non-identity morphisms. To describe the action of $\cW$  on the plane $\fh^*=\ttk^{1|1}$, take the objects to be $\pm a$ where $a \neq 0$ is a point in $\fh^*$.  The functor sends  $\pm a$ to the line $L$ they span, and the non-identity morphisms add or subtract $a$ to a point on $L$. The invariant ring $S(\fh)^{\cW}$ consists of all polynomial functions that are constant on $L$. So $S(\fh)^{\cW}= \ttk + TS(\fh)$ where $S(\fh) =\ttk[S,T]$ and $L$ is the line defined by $T$. The algebra $S(\fh)^{\cW}$ is a well-known example of a non-Noetherian domain of Krull dimension two.   
This answers a question of S. Paul Smith.} 
\eexa

\subsubsection{The Duflo-Serganova functor. } \label{wbn}
We recall the definition of the Duflo-Serganova functor $DS_x$,   in the case that $x$ is a root vector such that $[x,x]=0.$ 
For $M$ as above, set 
 $M_x =\Ker_M x/x M$.  Then 
$M_x$ is a module for 
$\fg_x$ where   
$\fg_x =\Ker \ad x/ \Im \ad x,$ which is a Lie superalgebra of smaller rank  than $\fg$. Although 
$\fg_x$ is a quotient of a subalgebra of  $\fg$, it was shown in 
\cite{DS} Lemma 6.3, that it embeds in $\fg$ at least in the KM case.  We follow the description of 
\cite{HR} 
 section 2C in the case that $x$ is a  root vector corresponding to the isotropic root $\gb$.

\bl \label{btal} 
Let $x$ be a root vector corresponding to the isotropic root $\gb$.
We can identify $\fg_x$ with a subalgebra of $\fg$ such that the centralizer $\fg^x$ of $x$ in $\fg$ is a semidirect sum $\fg^x= \fg_x \ltimes [x,\fg]$. The subalgebra $\fg_x$ is given by
\be \label{ktal}\fg_x= \fh(\gb)\op \bop_{\ga\in \Gd(\gb)}\fg^\ga,\ee where 
$\Gd(\gb) = \{\ga\in \Gd|(\ga,\gb) =0, \ga\neq \pm\gb\}$, 
and $\fh(\gb) = \span \{h_\ga|\ga \in \Gd(\gb)\}$.  We have $\fh(\gb) = \fg_x\cap \fh$.
\el
\noi 
\bl \label{ndg} The restriction of $(\;,\;)$ to $\span\Gd(\gb) $ is non-degenerate.
\el
\bpf Let $\{\ga_1,\ldots, \ga_r\}$ be a basis for $\span\Gd(\gb) $ and suppose  $\gs\in \fh^*$ with 
 $(\gs,\gb)=1$.  For $i\in [r],$ let  $\ga_i' = \ga_i - (\gs,\ga_i)$.  The Gram matrix of  $(\;,\;)$ using the basis 
$\{\ga_1,\ldots \ga_r,\gb,\gs\}$ for $\fh^*$ 
has the same determinant as the Gram matrix using the basis where each 
$\ga_i $ is replaced by $\ga_i'$.  Since the latter matrix is block diagonal and $(\ga_i, \ga_j) = (\ga_i', \ga_j')$, the result follows.
\epf \noi 
Under the isomorphism $\fh^* \lra \fh^*$,   $\span \Gd(\gb)$ maps to  $\span \{h_\ga|\ga \in \Gd(\gb)\}$ which we identify with $\fh(\gb)$. We have $\span \Gd(\gb) = \fh(\gb)^*$. 
If $(\gs(\gb),(\gb))=1$, we have an orthogonal direct sum 
\be \label{cut} 
\fh^*  = \fh(\gb)^* \op \span\{ \gb,  \gs(\gb)\}.
\ee
\noi In the table below, taken from \cite{GHSS} Section 4.5, the first column lists the Lie superalgebras we consider, $x$ is a root vector, and the second column describes $\fg_x$ up  to isomorphism. If $\fg=\fp(n)$, we assume that $x\in \fg^{\gep_1+\gep_2}$ in the notation of \cite{M}, 2.4.1.  
    \[ \begin{tabular}{|c|c|}\hline
$    \fg$ &$\fg_x$
\\ \hline\hline 
$    \fgl(m|n)$ &$\fgl(m-1|n-1)$ \\ \hline  
$    \fsl(m|n), \; m \neq n$ &$\fsl(m-1|n-1)$ \\ \hline 
$    \osp(m|2n) $ &$\osp(m-2|2n-2)$ \\ \hline
$    D(2|1;\ga)$ &$\ttk$\\ \hline 
$    F(4)$ &$\fsl(3)$ \\ \hline  
$G(3)$ & $\fsl(2)$ 
\\ \hline
$\fq(n)$&$\fq(n-2)$
\\ \hline   
$\fp(n)$&$\fp(n-2)$
\\ \hline 
    \end{tabular}  \]
 Since we use induction on  the rank as a proof technique, we make some brief remarks about the base cases.  These occur when $\fg_x$ 
has no isotropic roots.  Then the Weyl groupoid for $\fg_x$ is just the Weyl group. The desired results  hold  trivially in these cases. 
By convention $    \fgl(0|k) =    \fgl(k|0) =  \fgl(k)$, and other low dimensional entries in the column  for $\fg_x$ should be interpreted similarly. 
\\ \\
The main use we make of the functor 
$DS_x$ is to construct certain algebra  maps $ds_x$ in Subsections \ref{wxn} and \ref{wyn}. 
These algebra maps can be realized in two equivalent ways, and both realizations can be found in the literature, at least implicitly.  Both methods can be easily adapted  to obtain 
maps on suitable overrings as in the diagram preceding Lemma \ref{doe}. 
 In the situation of  Lemma \ref{btal}, $\fh^*$ has been  identified with a subspace of $\fh^*$, so we have a resriction map   $\res: S(\fh) \lra S(\fh_x), f \mapsto f|_{\fh^*}$. 
This map descends to $S(\fh)^W$  and  $I(\fh)$ and we have   $\res = ds_x: I(\fh) \lra I(\fh_x)$. On the other hand the algebra $I(\fh)$ can be defined in terms of polynomials satisfying certain partial evaluation conditions,  and then it is natural to view $ds_x$ as evaluation $\ev:I(\fh) \lra I(\fh_x)$.  Similar remarks apply to $ds_x: J(G) \lra J(G_x)$. 
\subsubsection{Harish-Chandra pairs and representations of Lie supergroups.} \label{whn}
An (algebraic) Lie supergroup $G$ can be defined using a Hopf superalgebra $\cO(G)$ and the functor of points. Denote  the
category of supercommutative $\ttk$-algebras ${\bf Alg}$.  Then for any supercommutative algebra $A$, the $A$-points of  $G$  are given by 
 $$G(A) = \mo _{\bf Alg
}(\mathcal{O}(G),A).$$ 
The Hopf superalgebra structure on $\cO(G)$ induces a natural group structure on $G(A)$. Another  useful approach (due to Kostant \cite{Ko} and  Koszul \cite{Kos}) to $G$ uses Harish-Chandra pairs. Following \cite{CCF} Definition 7.4.1, a {\it Harish-Chandra pair} (HCP) is a pair $(G_0,\fg)$
where consisting of an algebraic group $G_0$ and a Lie superalgebra $\fg$
such that  
\bi \itema $ \fg_0 = \Lie G_0,$ 
\itemb There is a representation $\gr$ of $G_0$ on $\fg$  such that $\gr(G_0)|_{\fg_0} =\Ad$ and the differential of $\gr$ acts on $\fg$ as the adjoint representation: for $X\in\fg_0, Y\in \fg$,
$$ d\gr(X)Y = [X,Y].$$
\ei 
We write $(G_0,\fg,\gr)$ if we want to stress the role of $\gr$. Harish-Chandra pairs form a category, for morphisms see \cite{CCF} Definition 7.4.2. 
Given a  Lie supergroup  $G$, we obtain a HCP $(G_0, \Lie G, \Ad)$.  Conversely given a HCP, we can put a Hopf superalgebra structure on
\be \label{hss} \cO(G):=\Hom_{U(\fg_0)}(U(\fg),\cO(G_0)),\ee
 and thus obtain a Lie supergroup $G$.  
In \eqref{hss}, $ U(\fg)$ is naturally a left $U(\fg_0)$-module, and the action of $x\in \fg_0$ on $\cO(G_0)$ is given by a left-invariant differential operator, \cite{CCF} page 125.  
This is the basis for the following result, see \cite{CCF} Theorem 7.4.5.
\bt \label{wln} There is an equivalence of categories between the category of Lie supergroups and  the category of  Harish-Chandra pairs.
\et\noi 
The connection between representations of  Harish-Chandra pairs and the corresponding Lie supergroup was explained in \cite{CC}. A {\it representation} of a HCP $(G_0,\fg)$ consists of the following data
\bi \itema A representation $\gS: G_0 \lra GL(V),$ 
\itemb A representation $\gs:\fg\lra \fgl(V)$,  such that $d\gS = \gs|_{\fg_0}$ and a further compatabilty condition given in \cite{CC} Definition 6 (2) holds.
\ei 
\bt \label{wmn} There is a bijection between representations of the HCP
$(G_0,\fg)$ 
and representations of $G$ on $V$. 
\et
\bpf See \cite{CC}  Proposition 6. \epf
Let $\cF_\fg$ (resp. $\cF_G$) 
 be the category of finite dimensional $\Z_2$-graded $\fg$-modules (resp. $G$-modules), and $\Pi$ the parity change functor on $\cF_\fg$. 
The extra condition mentioned above essentially says that the representations of $G$ that arise correspond to representations of $\fg_0$ that are integrable as  $G_0$-modules.  For example, suppose that $\fg_0$ has a non-trivial center.  If $\fa$ is an abelian Lie algebra and $T$ is a torus with $\Lie T = \fa$, then integrable representations of $\fa$ correspond to elements of the character group $X(T)$. 
Hence   $\fg_0$ has many more finite dimensional modules than 
$G_0$, and these modules can be induced to give   graded modules for $\fg$. Thus the category  $\cF_G$ is more manageable than $\cF_\fg$. In addition, the results we quote from \cite{HR}, \cite{IRS} and \cite{Re} work best for supergroups. 
Let $\cF=\cF_G$ or  $\cF_\fg$.  The Grothendieck group $K(\cF)$ of $\cF$ is the  free abelian group on the symbols $[M]$, where $M$ is an object in $\cF$ modulo the relations that $[M] = [L] +[N]$ whenever there is an exact sequence 
$0\drar L\drar M\drar N \drar 0$
in $\cF$.   Since $\cF$ is closed under tensor product, $K(\cF)$  has a ring structure determined by $[M ][N] =[M\ot N]$.  
The {\it supercharacter ring} (resp. {\it character ring}) 
 is the factor ring of $K(\cF_\fg)$ by the ideal generated by all $[\Pi(M)]+[M]$ (resp.  $[\Pi(M)]-[M]$) where $M\in \cF_\fg$. For the classical Lie superalgebras $\fg$ listed in Subsection \ref{Cls}, we define  $J(\fg)$ to be the   character ring  of  $\fg$ if $\fg$ has type $\fq$ and the  supercharacter ring in all other cases. The $\Z$-algebra $J(G)$  is the factor ring of $K(\cF_G)$ defined in a similar way.  The  Lie superalgebras $\fg$ listed in \ref{Cls} such that $\fg_0$ has a non-trivial center  are $\fg= \fgl(m|n),\fsl(m|n)$ and $\fg=\osp(2|2n).$  In these cases $J(\fg)$ is graded by an abelian group, and  $J(G)$ is the identity component of this grading, see \cite{SV2} Propositions 7.1, 7.2, 7.3 and 7.5. 
\\ \\
\noi  
We mention an immediate consequence of 
the weak  Nullstellensatz for $Z(\fg)$, Theorem \ref{wn}. 
Let $\fg$ be any  classical Lie superalgebra. A {\it central character } of $\fg$ is an algebra map $Z(\fg) \lra \ttk$, \cite{GHSS} Section 3.  Any central character $\chi$ determines, and 
is determined by the maxmal ideal $\Ker\chi$.  An element   $z\in Z(\fg)$ acts on any 
 Verma module $M(\gl)$ as a scalar $\chi_\gl (z)$ and  $\chi_\gl$ is called    the {\it central character  afforded} by $M(\gl)$, \cite{M} 8.2.4.  
The next result proves and extends 
Conjecture 13.5.1 from \cite{M}.
\bc \label{ccc} Any central character of $\fg$ is afforded by a Verma module. \ec
 \noi 
This paper replaces the unpublished preprint \cite{M19} which dealt only with the cases  
where $\fg=\fgl(m|n)$ or  
$    \osp(r|2n)$. 
I am grateful to    Maria Gorelik,  Hanspeter Kraft  and Shifra Reif for helpful correspondence.  

\section{General Results.} \label{se2}
Throughout  $[k]$ denotes  the set of the first $k$ positive integers and iff means if and only if. We start with some general results on commutative rings.
\noi  
\subsection{Certain Ring Extensions.} \label{cxe}
In this Subsection we assume the following hypothesis on the ring extension $A\subseteq B$.
\bh \label{xwl}   For some  
$T\in A$,  a non-zero divisor in $B$, we have  $TB \subseteq A$. \eh \noi 
If $I$ is an ideal of $A$, let $$\rad(I) = \{f\in A| f^n \in I, \mbox{ for some } n>0\}.$$
 If $I=\rad(I)$ we say $I$ is a {\it radical ideal}. 
Suppose $R$ is a  commutative ring  and $T\in  R$, is a non-zero divisor. Denote the localization by $R_T$. 
\bl  \label{owl}  Assume Hypothesis \ref{xwl}.
\bi \itema We have an equality of localizations $A_T = B_T.$
\itemb
If $P$ is a radical ideal
of $A$  with $T\in P$, 
then $TB \subseteq P$. \ei
\el
\bpf We leave (a) as an exercise. Then (b) holds since $(TB)^2
\subseteq A P = P$.
   \epf\noi 
If $L$ is an ideal of $R$ define the {\it  extension} of $L$ to $R_T $ to be $L^e = L_T$, the localization of $L$.   If $M$ is an ideal of $R_T$, set $M^c= M \cap R$, the {\it  contraction} of  $M$ to $R$. Then 
extension  and contraction provide order preserving bijections  
\be \label{yok} \{P \in \Spec   R|  T \notin P \} \leftrightarrow \Spec   R_T.\ee
\noi  
This applies to $R= A$ or $B$ in Lemma \ref{owl}, so   we have a bijection

\be \label{yak} \{P \in \Spec   B| T \notin P\} \leftrightarrow  \{P \in \Spec   A| T \notin P\}
,\ee given by $P \mapsto P\cap A$.    From now on we identify the two sets on either side of \eqref{yok}.
\bt \label{cow} Let $\phi: A\lra A/BT= C$ be the natural map.  Then we have a disjoint union
\be \label{gnu} \Spec   A = \Spec B_T \cup \phi^{-1} (\Spec   C).  \nn\ee
\et
\bpf As noted above we have $A_T =B_T $.    For a prime ideal $Q$ of $A$ there are two possibilities.  If $T\notin Q$, then $Q_T$ is a prime ideal of $A_T=B_T$ such that $Q_T\cap A=Q.$ If $T\in Q$, then $Q$ is the inverse image of the  prime ideal of $Q/BT$   under $\phi.$
\epf
\bc  \label{pig}  If         $m$ is a maximal ideal of  $A$ and $T\notin m$, the  maximal ideal $M$ of  $B$
given by $M= m_T\cap B$ satisfies  $m = M \cap A$.    Also $mB \neq B$.     \ec \bpf    This first statement follows from Equation  \eqref{yak}, and the second is an immediate consequence.         \epf
\noi In Corollary \ref{pig}, if $B=\cO(X)$ for an affine variety $X$, 
with $B$ a finitely generated $\ttk$-algebra, 
then  
  $m$ and $M$ have the same set of zeroes in $X$, so by the usual Nullstellensatz  $M  = \Rad mB$.  However $M$ can strictly contain $mB$, 
see  \cite{M22} Example 2.2.
\bl \label{hos} In the situation of Theorem \ref{cow}, suppose that the rings
$B_T $ and $ C$ satisfy  the ascending chain condition on radical ideals.  Then so does $A$. \el
\bpf    Let $R^{1} \subseteq R^{2} \subseteq \ldots $
be an ascending chain of radical ideals in $A$.
Then $R_T^{1} \subseteq R_T^{2}  \ldots $ is an ascending chain of radical ideals in $A_T=B_T$, so by assumption, there is an $m$ such that $R_T^{m} = R_T^{i}$ for all $i\ge m$.
Write $R_T^{m} = p_1 \cap p_2 \cap  \ldots  \cap p_r$ for some prime ideals of $B_T$.
  We can assume that this intersection is irredundant, and all the $p_k$ are minimal over $R_T^{m} $.
Then if  $P_k= p_k \cap A$ and $i\ge m$,  the $P_k$ are exactly the minimal primes over $R^{i} $
 that do not contain $T$. By \cite{E} Corollary 2.12 every radical ideal in a commutative ring is the intersection of (a possibly infinite number of)  prime ideals. So for $i\ge m$, write $R^{i}  $ as an intersection of prime ideals in $A$.  Say
\[ R^{i}  = P_1 \cap \ldots \cap P_r \cap  \bigcap_{j\in \Lambda_i} Q_j.
\]
 where the $Q_j$ are prime ideals that contain $TB$. Thus  $D_i =  \bigcap_{j\in \Lambda_i} Q_j$ is a radical ideal in $A$ containing $TB$, and if $\overline{D}_i = D_i/TB$, then $\overline{D}_{1} \subseteq \overline{D}_{2} \subseteq \ldots $ is an ascending chain of radical ideals of $C$.  So there is an $n \ge m$, such that $\overline{D}_{i} = \overline{D}_{n}  $ for all $i>n$ and  then $R^{i} = R^{n}  $. 
\epf
\bc \label{hxs}  With the same assumptions as the Lemma,
every radical ideal $I$ of $A$ is a finite intersection
of prime ideals. These prime ideals can be taken to be the prime ideals minimal over $I$.\ec  
\bpf See \cite {K} Theorem 87.\epf

\subsection{Towards the Weak Nullstellensatz.} \label{twn}
Consider a  ring  extension  $A \subseteq B$, and $X\subseteq \Max A$.  We say the {\it Relative Weak Nullstellensatz} (RWN) for $X$ and the pair $(A, B)$, if for every $m \in \Max A$ we have $mB \neq B$. If the RWN holds in the case $X=\Max A$, we say  simply that 
 RWN  holds for the pair $(A, B)$.  
If $A, B$ are $\ttk$-algebras, with $B$ finitely generated  and $M \in \Max B $, then $B/M \cong \ttk$, so if  $m=M\cap A$, then \be 
\label{sgs4}  \ttk \hookrightarrow A/m \cong (A +M)/M \subseteq B/M \cong\ttk. \ee  Thus $m \in \Max  A$ and we have a map 
 \be \label{sgs51} \psi:\Max B \lra  \Max A, \quad M \mapsto M\cap A . \ee  
If RWN  holds for the pair $(A, B)$, and $m \in \Max A$, then  by Zorn's Lemma  there is  a $M \in \Max B $ containing $mB$, we have $M\cap A = m$, so the map in \eqref{sgs51} 
is surjective.\\ \\
Consider the following diagram of ring homomorphisms where the vertical maps are inclusions, and 
${\phi}$ is the restriction of ${\Phi}$ to $A$.

\begin{xy}
(0,20)*+{}="f";
(40,20)*+{A}="a";
 (100,20)*+{A'}="b"; (40,0)*+{B}="c";
(100,0)*+{B'}="d"; {\ar "a";"b"}?*!/_2mm/{\phi}; {\ar "a";"c"};
{\ar "b";"d"};
{\ar@{>} "c";"d"};?*!/_2mm/{\Phi};
 \end{xy}
\bl \label{doe} Suppose $\Phi$ and $\phi$ are onto and  {\rm RWN} holds 
for the pair $(A', B')$. Then the RWN holds for the pair $(A, B)$  
and the set 
$X= \phi^{-1}(\Max A)$.\el
\bpf 
Suppose 
$m \in  X $, that is $m$ is a maximal ideal of $A$ containing $\Ker \phi$. 
Then 
 $\phi (m)$ is  a maximal ideal of $A' $, so by assumption  $J=\phi(m)B'$ is a proper ideal of $B' $. Since $\Phi(mB )=J,$ it follows that $mB $
is a proper ideal of  $B$.
\epf
\bl \label{doy} Suppose the finite group  $G$  acts on $S$ and  assume $|G|$ is a unit in $S$.  If {\rm RWN} holds 
for the pair $(R, S)$, then the RWN holds for the pair $(R^G, S^G)$.\el
\bpf In this situation we have a Reynolds operator $R\lra  R^G$.  Regarded as a $G$-module $R^G$ is the isotypic component of $ R$ corresponding to the trivial module.  Let $R^+$ be the sum of the remaining isotypic components.  Thus 
 $ R=R^G\op  R^+$ 
is a decomposition into $R^G$-modules. If $m \in 
\Max  R^G$, then $mR= m \op mR^+$ is a proper ideal of $R$, so  by assumption  $mS$ is a proper ideal of $S$. The result follows since  $ m S^G= S^G$ implies   $mS=S $.
\epf
\noi RWN does not hold always when $S =\ttk[x,y]$ is a polynomial algebra in 2 variables: if $R = \ttk[x] + (xy-1)S$, then $m=x\ttk[x] + (xy-1)S$ is a maximal ideal of $R$ such that $mS=S$.
\section{Applications.} \label{se3}
In this Section  if $\A=I(\fh)$ or   $J(G) \ot_{\Z}\ttk$ we consider  a ring extension $\A\subset \B$ and an element $T\in \A$ such that $T\B \subset \A$. Thus Hypothesis \ref{xwl} holds. Because we will need to change the Lie superalgebra, we sometimes   write $\A=\A(\fg)$. 
 The algebra $\B$ is either $S(\fh)^W$ or $\cO(\mathbb{T})^W$,
so is in particular a finitely generated  $\ttk$-algebra.  To prove the weak Nullstellensatz we consider a map  with kernel $T\B$, 
\be \label{jtal} \A(\fg)\lra  \A(\fg_x),\nn\ee
where $\fg_x$ is as in the table from Subsection \ref{wbn}, and use induction on the rank of $\fg$. 

\subsection{Application to $I(\fh)$.} \label{wxn}
By \cite{GHSS} Section 3, the functor   
$DS_x$  induces an algebra map $Z(\fg) \lra Z(\fg_x)$. Using the Harish-Chandra map, this gives rise to a map     $ds_x:I(\fh) \lra I(\fh_x)$.  Because of \eqref{ktal} and the above remarks, this map can be thought of as restriction of functions
$\res(f) = f|_{\fh_x^*}$. 
\subsubsection{The Kac-Moody  case.}\label{KMc}
In the KM case, if $\fg_1$ is a simple $\fg_0$-module, let $\gO=\Gd_{iso}$.  Otherwise 
$\fg_1= \fg_1^+ \op \fg_1^-$ is a direct sum of two simple $\fg_0$-modules, 
 and  we 
let $\gO=\Gd_{1}^+$ be the set of roots of $\fg_1^+$.  Then $W$ acts transitively on $ \gO$.  Define 
 $T=\prod_{\beta\in\gO} h_{\beta}$. Then  $T$ is $W$-invariant. 
 We describe  an automorphism  $\gs$ on  $\fg_x$ 
which we are
  going to use below. 
\bi
 \itema
 For $\fg=\osp(2m|2n)$,  then  $\gs$ is the diagram  automorphism of $\fg_x$ described in \cite{M} Lemma 5.5.12. This automorphism preserves $\fh$. 
 \itemb
 For $\fg=F(4)$, then 
  $\gs$ is induced by an involution of the Dynkin diagram for $\fg_x=\fsl(3)$.
\itemc
  for $\fg=D(2|1;\ga)$ we have $\fg_x=\ttk$ and  $\gs=-\Id$.
  \itemd In all other cases  $\gs=\Id$.
  \ei
\bt 
\label{thmtheta} We have 
\bi 
 \itema
 $T\cS(\fh)^W\subset I(\fh)$ and the
kernel of $\res: I(\fh)\to I(\fh_x)$ equals 
$T\cS(\fh)^W$.
 \itemb
We have $\res(I(\fh)= I(\fh_x)^{\gs_x}$. \ei
\et
\bpf
By construction $T$ is $W$-invariant and $T$ vanishes on all hyperplanes $\Pi_{\ga}$. Thus by \eqref{ltal}, 
$T\cS(\fh)^W\subset \Ker\res$. Suppose $x =e _\gb$ is a root vector and  $f \in\Ker\res$. 
Then $h _\gb$ divides $f$ in $\cS(\fh)$, so by  $W$-invariance, so does 
$T$.  Write $f=gT$ with
$g\in \cS(\fh) $.  Since $f,T$ are  $W$-invariant,  so is $g$
proving (a).  For (b) see \cite{Gor1} Theorem 6.4.
\epf 
 \subsubsection{Algebras of type $\fq$.}\label{pqc}
There are four related Lie superalgebras of type  $\fq$. We refer to \cite{Re} 2.1, \cite{CW}  for full details, but mention that $\fq(n)$ is constructed as a subalgebra of $\fgl(n|n)$, and has derived 
subalgebra $\sfq(n)$. We have $\fq(n)_0 \cong \fgl(n)$ and  $\sfq(n)_0\cong \fsl(n)$.  The algebras $\fq(n)$,  $\sfq(n)$ have factor algebras 
 $\pq(n)$
and $\psq(n)$ respectively.  The algebra $\psq(n)$ is simple for $n\ge 3$. 
By  
\cite{Gq} Theorem 13.1.
\be \label{wpn1}
Z(\fq(n))\cong \{ f \in S(\fh)^W|  f|_{x_n =- x_{n-1}=t} \mbox{ is independent of } t\} := I_{n}.
 \ee Now \eqref{wpn1}   shows that if $T=\prod_{i<j} (x_i + x_j)$, then $T\in \A = I_n$. Set  $\B =  S(\fh)^W$.  Define 
f $\ev: I_{n}\to I_{n-2}$ by $\ev(f) = f|_{x_n =- x_{n-1}=0}.$ 
\bt 
\label{thm2} We have 
\bi 
 \itema
 $T\B\subset I_{n} $ and the
kernel of $\ev: I_{n}\to I_{n-2}$ equals 
$T\B$.
 \itemb
We have $\ev(I_{n}) =I_{n-2}$ \ei
\et
\bpf (a) follows as in  the proof of Theorem \ref{thmtheta}  (a).  For (b) see \cite{Gor1}  Proposition 5.8.3.
\epf \noi To unify notation, we set $I(\fh) =I_{n}$ if $\fg=\fq(n)$. In this case, by
 \cite{M} Lemma 12.1.1.
\be \label{etal}  I(\fh)= \{f\in S(\mathfrak{h})^W| f(\lambda) = f(\lambda + t \alpha) \mbox{ for all }\ga \in \Gd^+, \gl \in \Pi_{\ga} \mbox{ and } t \in \ttk\}\ee
 \noi 
 \subsubsection{Algebras of type $\fp$.}\label{pc}
We are interested in 
two Lie superalgebras of type  $\fp$. We refer to \cite{CW} 
 for details, but recall  that $\fp(n)$ is constructed as a subalgebra of $\fgl(n|n)$, and has derived 
subalgebra $\sfp(n)$. We have $\fp(n)_0 \cong \fgl(n)$ and  $\sfp(n)_0\cong \fsl(n)$.  
The algebra $\psq(n)$ is simple for $n\ge 3$. 
If  $\fg =\fp(n)$, then by \cite{Gor} Theorem 4.1, $ Z(\fg)$ contains 
an element $t$ such that  $t^2=0$ and there is a linear isomorphism  $Z(\fg) \cong  \ttk \op tS(\mathfrak{h})^W$.  This means that $Z(\fg)$ has a unique prime ideal and the methods of Subsection \ref{cxe} do not apply to the study of $Z(\fg)$.  For this reason $\fp(n)$ is  excluded in some results. Corollary \ref{ccc} holds when  $\fg=\fp(n)$, since there is only one central character.  The ring 
 $J(G)$ when $G$ is the Lie supergroup $P(n)$ remains of interest.

\subsection{Application to $J(G) \ot_{\Z}\ttk$.} \label{wyn}\noi 
\subsubsection{Lie supergroups}
 We make a  basic hypothesis, compare \cite{GHSS} Section 8. 
\bh \label{bns}{\rm Assume one of the following holds
\bi \itemi If  $ \fg = \fgl(m|n) $, set $G_0=GL(m)\ti GL(n) $.%
\itemii If  $ \fg = \osp(r|2n)$, set  $G_0=SO (r)\ti SP(2n) $, where $r=2m+1$ or $2m$.
\\
  In  cases (i) and (ii) set 
\be\label{tej}P= \sum_{i=1}^m \Z\gep_i \op \sum_{i=1}^n \Z\gd_i.\nn \ee 
\itemiii If  $ \fg = \fsl(m|n) $ $m\neq n$, set $G_0=\{(A,B)\in GL(m)\ti Gl(n)| \det A = \det B \} $ and 
let $P$ be the root lattice of $\fg$.
\itemiv If  $ \fg = \fp( n) $, set $G_0=\{(A,A^{-t})\in GL(n)\ti GL(n)\} $.
\itemv If  $ \fg = \fq(n) $, set $G_0=\{(A,A)\in GL(n)\ti GL(n) \} $.
\ei}
\eh \noi In (iv) $A^{-t}$ is the transpose of $A^{-1}$. Suppose  (iv) or  (v) holds.  The    simple roots are 
$\ga_i =\gep_i - \gep_{i+1}, i \in [n]$.  We use the bilinear form on 
$ \fh^*$ given by 
$(\gep_i,  \gep_{j}) =\gd_{i,j}$ and set $P = \sum_{i=1}^{n}\Z\gep_i$. 
Recall that  $\mathbb{T}$ is a maximal torus in $G_0$.
\bl \label{abc}\bi \itemo\itema  
There is a
Harish-Chandra  pair $(G_0,\fg)$ and the Lie supergroup determined 
by \eqref{hss} is 
 $G=GL(m|n),  OSP(r|2n),  SL(m|n), P(n) $  or $Q(n)$ in  cases $($i$)$-$($v$)$ respectively.
\itemb There is an isomorphism  $P \lra X(\mathbb{T})$, which we write as $\ga\lra \tte^{\ga}$. 
Here $\tte$ is a formal symbol used to convert to multiplicative notation.
\ei 
\el \bpf The possible weights of 
$G_0$-modules (and hence also $G$-modules) are determined by $P$. Thus the result follows from the explicit description of the characters or supercharacters of these modules given in \cite{SV2}, \cite{IRS} and \cite{Re}.
\epf \noi
If $ \fg = \fp( n) $ or   $  \fq(n) $, and  $\ga=\gep_i - \gep_{j} \in \Gd_0^+,$ the set of  positive even roots, set  $\overline{\ga}=\gep_i + \gep_{j}$ and 
define in place of \eqref{rtal},
\be \label{rtak}  \Pi_{\ga}=
 \{ \lambda \in   \fh^* | (\lambda,\overline{\ga}) =0\}.
\ee  In the KM case, an {\it iso-set} is a linearly independent set of mutually orthogonal isotropic roots. If $ \fg = \fp( n) $ or   $  \fq(n) $, an  {\it iso-set} is  a linearly independent subset $A$ of $\Gd_0^+ $ such that $(\ga,\overline{\gb}) =0$ for all $\ga, \gb\in A$.
The non-standard convention that $\Gd_{iso} = \Gd_0^+ $ in types p and q, allows us to give unified statements of certain results.
\noi 
Let $P\subset \fh^*$ be as in Hypothesis \ref{bns}. For all $\ga \in P$, define the subtorus $\ttT_\ga $ 
of  $\mathbb{T}$ by 
\be \label{rtaj} 
\ttT_\ga = \Ker \tte^{\ga} \mbox{ in cases (i)-(iii), or }   \ttT_\ga = \Ker \tte^{\overline{\ga}} \mbox{ in cases (iv)-(v).} \ee
Composition gives a non-degenerate pairing 
$X(\mathbb{T}) \ti Y(\mathbb{T}) \lra \Hom(\ttk^*,\ttk^*) \cong \Z.$  The next result relates this pairing to the bilinear form on $P$.
 \bl \label{cwg3}
For each  root $\gb$ 
there is a unique $c_\gb \in Y(\mathbb{T})$ such that  
\be \label{cwg4} \tte^{\ga}c_\gb(t)= t^{(\ga,\gb)} \mbox{ for all roots } \ga 
\mbox{ and } t \in \ttk^*.\ee
\el 
\bpf First suppose   $G=GL(m|n) $ or
$    OSP(r|2n), $ with $r = 2m+1$ or $2m$ and 
let $p=m+n$. We show \eqref{cwg4} holds for all $\gb\in P$. 
It is convenient to relabel the ordered basis of $P$ as follows 
$$(\ga_1,\ldots,\ga_p):= (\gep_1,\ldots,\gep_m, \gd_1,\ldots,\gd_n).$$ Then the  dual basis is 
$$(\varpi_1,\ldots,\varpi_p):= (\gep_1,\ldots,\gep_m, -\gd_1,\ldots,-\gd_n).$$
There are unique $c_1,\ldots,c_p  \in Y(T)$ such that $\tte^{\ga_i}c_j(t)= t^{\gd_{i,j}}$. 
 If $ \gb \in P$, write $\gb= \sum_{i=1}^p a_j 
\varpi_j$ with $a_j\in\Z$.  Then set $ c_\gb = \prod_{i=1}^p c_j ^{a_j}$. It is enough to show \eqref{cwg4} for $\ga=\ga_i$, but in this case we have, for $t \in \ttk^*$ we have 
$$\tte^{\ga_i}c_\gb(t)= t^{a_i} =t^{(\ga_i,\gb)}.$$   
If   $G=SL(m|n) $ the distinguished set of simple roots is 
\begin{eqnarray*}
\alpha_i & = &  \epsilon_i - \epsilon_{i + 1}, \quad\quad i \in [m - 1]\\
\alpha_{m} & = & \epsilon_{m}  - \delta_{1}, \quad\quad   \\
\alpha_{m+i} & = & \delta_i - \delta_{i+1}, \quad\quad   i \in   [m - 1]
\end{eqnarray*}
Set $x_{i} = \tte^{\epsilon_i}, y_{i} = \tte^{\delta_i}$. Then for a simple root $\ga$, 
$\tte^\ga =x_{i}x_{i+1}^{-1}, x_{m}y_{ 1}^{-1}$
and $y_{i}y_{i+1}^{-1}$ in the three cases listed above. For a root $\gb$ define 
$c_\gb \in Y(\mathbb{T})$, by
$$c_\gb(t) = (t^{(\gb, \gep_1)}, \ldots , t^{(\gb, \gep_m)}, t^{(\gb, \gd_1)}
, \ldots ,  t^{(\gb, \gd_n)}).$$ It suffices to show \eqref{cwg4} for $\ga, \gb$ for simple roots and this is done with a short computation. 
 \epf   \noi
 \bl \label{cwg31} If 
$G=P(n) $  or $Q(n)$, then for all  $\gb\in \Gd_0^+$ 
there is a unique $c_\gb \in Y(\mathbb{T})$ such that  
\be \label{cwg41} \tte^{\overline{\ga}}c_\gb(t)= t^{(\overline{\ga},\gb)} \mbox{ for all   } \ga \in \Gd_0^+
\mbox{ and } t \in \ttk^*.\ee
\el 
\bpf Define $c_\gb \in Y(\cT)$ by $c_\gb(t)=(1, \ldots, t, \ldots , t^{-1}, \ldots ,1)  $, where $t, t^{-1}$ are in positions $i$ and $j$.  The result holds by a short computation.
\epf
In the KM case, it was shown by Hoyt and Reif, \cite{HR} Proposition 8, that the functor  
$DS_x$ induces a ring homomorphism $ds_x:J(G)\lra J(G_x).$ 
Our next aim is to describe the image and kernel of the map $ds_x$, and deduce that the inclusion 
$\A=J(G) \ot_{\Z}\ttk \subseteq \B=\Z[X(\mathbb{T})]^W \ot_{\Z}\ttk
$ satisfy Hypothesis \ref{xwl}.  Though not strictly necessary for all the results we cite, we extend scalars to $\ttk$. This brings out the analogy with $I(\fh)$. Note that 
$\Z[ X(\cT)]\ot_\Z \ttk=\cO(\cT).$
\subsubsection{Properties of  $ds_x$.} \label{wcn} In the the KM case, set 
${\gr_\iso}=\frac{1}{2}\sum_{\ga \in \Gd^+_\iso} \ga$. 
\bt \label{wkn}\bi \itema There is an element $T\in\B$ such that 
$\Ker ds_x =T\B\subset \A$. In particular  Hypothesis \ref{xwl} holds.
\itemb Let $G$ be one of the Lie supergroups $SL(m | n), m \neq n, GL(m | n), $ or 
$OSP(r| 2n)$. Then $ds_x: J(G)\lra J(G_x)$ is surjective.\ei\et 
\bpf (a)  Define $T$ by 
$$T = \tte^{\gr_\iso}\prod_{\ga\in \Gd^+_\iso}(1-\tte^{-\ga}).$$  In the notation of \cite{HR} Lemma 16  we have $T = k({\gr_\iso})$.  So the result follows from the cited Lemma and \cite{HR}  Theorem 17. Also (b) 
 is \cite{HR}  Theorem 20.
\epf  \subsubsection{The case of  $Q(n)$.} \label{wfn}
The case   $G=Q(n)$ is studied in \cite{Re}. 
Here are some relevant results. In this and the following Subsection, let  $\mathbb{T}_n$ be a maximal torus in $Q(n)_0=GL(n)$. Then  \cite{Re} 
Proposition 3, gives an explicit description of a ring $J_{n}$  isomorphic to  $J(Q(n))$.  For our purposes it is convenient to extend scalars to $\Z[\frac{1}{2}]$. Thus we define 
\be \label{wpn}
J_n = \{ f \in \Z[\frac{1}{2}][X(\mathbb{T}_n)]^W|  f|_{x_n =- x_{n-1}=t} \mbox{ is independent of } t\}.
 \ee With this definition $J_n \cong J(Q(n)) \ot_\Z \Z[\frac{1}{2}].$
Using this isomorphism, 
$ds_x$ identifies with  an evaluation map $\A = J_{n}\lra J_{n-2}=\A'$, given by $\ev(f) = f|_{x_n = -x_{n-1}=1}$.
Note that  $\A$ is contained in $\B = \Z[\frac{1}{2}][X(\mathbb{T}_n)]^W$. 
So we have a diagram as in Subsection \ref{twn}, where the map  $\phi=\ev$, and the map $\Phi:
\B \lra \B' = \Z[\frac{1}{2}][X(\mathbb{T}_{n-2})]^{W'}
$ is defined by 
\be \label{wrn} \Phi(f) = f|_{{x_n} = - x_{n-1} =1}.\nn\ee
\bt \label{wjn9}
\bi \itema The map $ \ev:J_{n}\lra J_{n-2}$ is surjective. 
\itemb Define 
$T \in \B=\Z[\frac{1}{2}][X(\mathbb{T}_n)]^W$, by $T= \prod_{i<j} (x_i + x_j)$.  Then
 $ T\B \subset \A$. That is Hypothesis \ref{xwl} holds. Furthermore $\Ker \ev = T\B.$
\ei 
\et
\bpf For (a) see \cite{Re} Proposition 7.  Since clearly $T$ is fixed by $W$ and $ \ev(T)=0,$ we have
$ T\B \subset \A$ and $ T\B \subseteq \Ker \ev.$
The other inclusion in (b) follows from the proof of \cite{Re} Proposition 6. 
Without  extending  scalars, we would not have $ T\B \subset \A$.
\epf
\subsubsection{The case of  $P(n)$.} \label{wgn}
The case of  $G=P(n)$ is studied in \cite{IRS}. We define 

\be \label{wqn}
J_n = \{ f \in \Z [X(\mathbb{T}_n)]^W|  f|_{x_n =x_{n-1}^{-1}=t} \mbox{ is independent of } t\}.
 \ee By  \cite{IRS} Theorem 1.0.1 $J_n \cong J(P(n)),$ the super character ring of $P(n)$. As before this allows us to identify 
$ds_x$  with  an evaluation map $\A = J_{n}\lra J_{n-2}=\A'$, given by $\ev(f) = f|_{x_n = x_{n-1}^{-1}=1}$.
Note that  $\A$ is contained in $\B = \Z [X(T_n)]^W$. So we have a diagram as in Subsection \ref{twn}, where the map  $\phi=\ev$, and the map $\Phi: 
\B \lra \B' = \Z[X(\mathbb{T}_{n-2})]^{W'}
$ is defined 
 by $\Phi(f) = f|_{x_n = x_{n-1}^{-1}=1}.$
\bt \label{wjn}
\bi \itema The map $ \ev:J_{n}\lra J_{n-2}$ is surjective. 
\itemb Define $T \in \B=\Z [X(\mathbb{T}_n)]^W$, by $T= \prod_{i<j} (1-x_i x_j)$.  Then
such that $ T\B \subset \A$. That is Hypothesis \ref{xwl} holds. Furthermore $\Ker \ev = T\B.$
\ei 
\et
\bpf For (a) see \cite{IRS} Theorem 4.0.1.   Statement (b) follows from the proof of \cite{IRS} Proposition 3.2.1.
\epf

\section{The  Nullstellensatz.}\label{sn} \noi 
In this Section we assume one of the following holds, compare Notation \ref{bn}
\bi \itemi $\fg\neq\fp(n)$, $\A= I(\fh), \; 
 X = \fh^* \mbox{ and } \B= \cO(X)^W =S(\fh)^W.$
\itemii $\A= J(G) \ot_{\Z}\ttk, \;  X = \mathbb{T}, 
\B= \cO(X)^W$ and $G$ is not exceptional KM. 
\ei\noi 
\subsection{Maximal Ideals.} \label{wn9}
\noi First we prove the analog of the weak Nullstellensatz.
\bt\label{hen}
  If $m$ is a maximal ideal of  $\A$, then there is a maximal ideal $M$ of  $\B$ such that $m = M \cap \A$.\et
\bpf We give the details in case (i).   
If $T\notin m$ this follows from Corollary \ref{pig}. Assume is KM   and let ${W_x}$ be the Weyl group of ${\fg_x}$. 
By induction the conclusion holds with $(\A,\B)$ replaced by the pair $(I(\fh_x), S(\fh_x)^{W_x}).$  Thus by 
Lemma \ref{doy} it holds also for  $(I(\fh_x)^{\gs}, S(\fh_x)^{W_x,\gs})$ where $\gs=\gs_x$ is as in Theorem \ref{thmtheta}.
Thus if  $T\in m$, the result holds by Lemma \ref{doe}, taking  $M$ to be any maximal ideal containing $m\B$. The same argument works for  $\fg=\fq(n)$ using Theorem \ref{thm2}.  The proof in case (ii) is similar.
\epf

\bt \label{wn} If  $\gl \in X$, set 
$$m_\gl = \{ f\in \cA| f(\gl)=0\}.$$
 If  $m$ is a maximal ideal of $\A$, then $m = m_\gl$ for some   $\gl\in X$.
\et
\bpf This follows from Theorem \ref{hen}.   Note that  maximal ideals in $\B$ correspond to $W$-orbits in $X$. 
\epf
\noi For $\gl \in X,$ let $M_\gl\in \Max \cO(X)$ be the  ideal of functions  vanishing at $\gl$.
\noi \bc   \label{csn} We have \bi \itema 
 There are surjective maps $$\Max\cO(X) \lra \Max \B  \lra \Max \A.$$   
\itemb Denote the composite of the maps in $(a)$ by $\gs$ and suppose that $\sim$ is
the smallest equivalence relation on $X$ such that
\begin{itemize}\itemi $w\lambda \sim \lambda$ for all $w \in
W.$
 \itemii If $X=\fh^*$ and  $(\lambda , \alpha) = 0$ for $\ga\in \Gd_{iso}$, then $\lambda \sim \lambda + t\alpha$ for all $t \in \ttk.$
\itemiii If $X=\cT$ and  $\lambda \in \ttT_\alpha $ for $\ga\in \Gd_{iso}$,  then $\lambda \sim 
c_\alpha(t)
\lambda  $ for all $t \in \ttk.$
\end{itemize}\vspace{0.2cm}
Then $ M_{\mu } \in\gs^{-1}(m_\gl)$ iff $\lambda \sim \mu.$ \ei\ec
\bpf In (a) the map $\Max \B\lra \Max \A$ is  given by $M\lra  M \cap \A$, and the other map is defined similarly.  The first map is surjective since $
\B$ is the fixed ring of  $\cO(X)$ under the action of the finite group $W$.  The second map is surjective  by Theorem \ref{hen}. If $ X = \fh^* $, part (b) is a reformulation of \cite{M} Theorem 13.5.4, and the proof when $ X = \cT $ is similar.
\epf \noi 
\subsection{The Strong Nullstellensatz.}
We deduce the strong Nullstellensatz from Theorem \ref{hen}. There is a small problem because we need the analogous result when the pair $(\A, \B) $ is replaced by $(\A \ot_\ttk {\ttk[z]}, \B\ot_\ttk {\ttk[z]})$. This is easily taken care of by noting that the description of 
$\A$  given by  \eqref{ltal} and \eqref{yta} is independent of the algebraically closed field, and hence holds also for any field extension $K$ of $\ttk$.
\bc \label{dox} If  $m\in \Max \A \ot_\ttk {\ttk[z]}$, then $m\B\ot_\ttk {\ttk[z]}$  is a proper ideal of $\B\ot_\ttk {\ttk[z]}$. \ec
\bpf Let $S = {\ttk[z]} \bsk 0$.  If $m \cap S= \emptyset$, then $m_S \in \Max \A \ot_\ttk  {K}$,   
where 
$K={\ttk(z)}$, and the result   follows by the above remarks.  Otherwise, we have $z-\gl\in m$ for some $\gl\in \ttk$.  Now if $C= \A$ or  $\B$, set $C_\gl = C\ot_\ttk {\ttk[z]}/(z-\gl).$ Then $\overline{m}:= m/(z-\gl)   \in \Max \A_\gl .$  Now 
$$ \A \cong \A_\gl \subseteq \B_\gl \cong \B,$$  so by Theorem \ref{hen}  
$\overline{m} \B_\gl$ is a proper ideal of $\B_\gl$.
\epf

\noi If $I$ is a subset of $\A$, let
$\cV(I) = \{x\in X| f(x) =0 \mbox{ for all } f \in I\}$. 
Such a set is called an {\it superalgebraic set}. If instead
$I$ is a subset of $\B$, we say that $\cV(I)$ is an {\it algebraic set}
(or {\it closed}) in $X$. 
Thus any superalgebraic set is algebraic.
In addition if
$V$ is a subset of $X$, 
 set $$\cI_\cA(V) = \{f\in \cA| f(x) =0 \mbox{ for all } x \in V\}.$$ We will also need
$$\cI(V) = \{f\in \cO(X) | f(x) =0 \mbox{ for all } x \in V\}.$$
\bt \label{boa}  The maps $\cI_\cA$ and $\cV$ are inverse bijections between the set of superalgebraic sets in $X$, and the set of radical ideals in $\A$. Both maps are order reversing. \et

\bpf The key point is that $\cI_\cA(\cV(I)) \subseteq \rad(I)$. To show  $I$ is not  assumed finitely generated, we  repeat the well-known  ``Rabinowitsch trick",  
\cite{E} Theorem 4.19,
\cite{F} Chapter 1. Equation \eqref{bxa} below uses only finitely many elements from $I.$
Suppose $G\in \cI_\cA(\cV(I))$, and let $J$ be the ideal of $\A  \ot \ttk [z]$ generated by $I$ and $zG-1$. Then $\cV(J)$ is empty, since $G$ vanishes wherever all polynomials in $I$ vanish. Therefore by Corollary \ref{dox}, $1\in J$, and we can write
\be \label{bxa} 1 = \sum_{i=1}^r A_i F_i + B(zG-1),\ee 
where the $F_i $ are in $I$ and $B, A_i \in \cA  \ot \ttk [z].$  
Now set $Y=1/z$ and multiply by a large power of $Y$ to obtain
\[Y^N = \sum_{i=1}^r C_i F_i + D(G-Y)\] where $D, C_i \in \cA \ot \ttk [Y].$
The result follows by setting $Y=G$.
\epf
\bp \label{nit}
The maps $\cV,$ and $\cI_\cA$ satisfy  the following properties. Suppose that $E_\gl, V_\gl$ are subsets of  $\A$, and
$X$ respectively and that
$\fa, \fb$ are ideals of $\cA.$
\bi
\itema $\cV(\cup_{\gl\in \Gl} E_\gl) = \cap_{\gl\in \Gl}\cV( E_\gl)$.
\itemb $\cV(\fa \cap \fb) = \cV(\fa \fb) = \cV(\fa) \cup \cV(\fb)$.
\itemc $\cI_\cA(\cup_{\gl\in \Gl}V_\gl) = \cap_{\gl\in \Gl}\cI_\cA( V_\gl)$.
\ei\ep
\bpf Left to the reader. \epf

\subsection{Prime ideals and irreducible components.}
We say that  a superalgebraic set 
 is {\it irreducible}
 if it cannot be written as the union of two proper superalgebraic sets.
\bp  \label{hog}  If $I$ is a radical ideal of $\A$ and $V = \cV(I)$ is the corresponding superalgebraic set, then $I$ is prime iff $V$ is irreducible as a superalgebraic set.\ep
\bpf  The proof is completely analagous to the classical case \cite{F}
Chapter 1.5, Proposition 1, page 15.\epf
\noi
In general if
$I$ is a radical ideal of $\cA$, then using Corollary \ref{hxs}, we can write $I$ uniquely in the form $I =  P_1 \cap \ldots \cap P_r$,
where the $P_i$ are the prime ideals of $\cA$ which are minimal over $I$.  
Then $\cV(I) =  \bcu_{i =1 }^{r}  \cV(P_i)$ by Proposition \ref{nit}.
We call the superalgebraic  sets $\cV(P_i)$ the  {\it irreducible superalgebraic  components of $\cV(I)$.}  They are irreducible.

\bc Every superalgebraic set is uniquely a finite union of  irreducible superalgebraic components.
\ec

\bpf This follows  from  the Nullstellensatz and Proposition \ref{nit} (b).
\epf

\br {\rm Our approach to the Nullstellensatz was to relate ideals of $\A$ and $\B = \cO(X)^W =\cO(X/W).$ Therefore it might seem more naural to relate radical ideals in  $\A$  to their zero loci in $X/W$. However the connection between closed sets in $X$ and $X/W$ is well understood, when $\cO(X)$ is finitely generated,  \cite{DK} Section 2.3.1.  
If $U$ is a closed $W$-invariant subset of $X$, then $  \cI(U)$ is a  $W$-stable  ideal  of $\cO(X)$. Let $\pi:X\lra X/W$ be the natural map, corresponding to the inclusion of rings $\cO(X)^W \subseteq \cO(X)$. 
If $V $ is  closed  in $ X/W$, and $U$ is an irreducible component of $\pi^{-1}(V),$ then so is $wU$ for $w \in W$, and 
$\pi(U) = \pi (wU)$. In Section \ref{csc} it will be more convenient to work with closed sets in $X$, rather than  $X/W$. Several  closed subsets of $X $  we consider are identified under $\pi $, see Lemma \ref{AL}. 
}\er

\section{Weyl Groupoids and their Actions.} \label{sas}
\subsection{Actions of groupoids on varieties.}\label{awr} A groupoid $\fG$  can be
defined as a small category with all morphisms  invertible.
We denote the set of objects by $\mathfrak B$ which we call the {\it base}.
As in \cite{SV2}  we use the same notation $\fG$ for the set of morphisms  as for the groupoid itself.
For a variety  $X$ consider  the  groupoid 
$\mathfrak A(X)$ with  base all  subvarieties of $X$, and  morphisms all 
isomorphisms between subvarieties. We say the groupoid $\fG$ acts on   $X$  
if there is a functor 
\be \label{bwr}F:\fG \lra \mathfrak A(X).\nn\ee
If $g:b\lra b'$ is a morphism in $\fG$, set $s(g) = b$. 
Then for $\gl\in X,$ set $\fB_\gl =\{b\in \fB| \gl\in F(b) \}$ and $\fG_{(\gl)} =\{g\in \fG|s(g) \in \fB_\gl\}.$ If $g\in \fG_{(\gl)}$, we say $g$ {\it is defined at} $\gl $ and often write $g\gl=F(g)\gl$.
Set $\fG_{(0)}\gl = \{\gl\}$ and for $i \ge 0$, define 
$\fG_{(i+1)}\gl  = \{g\mu|\mu\in \fG_{(i)}\gl, 
g\in \fG_{(\mu)}
\}$.  The {\it $\fG$-orbit} of $\gl$  is 
$\fG \gl =\bcu_{i\ge 0} \fG_{(i)}\gl$. Informally, $\fG_{(i)}\gl$ is the set of points in $X$ we can reach from $ \gl$ using $i$ morphisms from $\fG.$ 
If $X$ is affine, the {\it invariant ring} for the action of $\fG$ on $X$ is 
 \[ \cO(X)^\fG =\{f \in \cO(X)| f \mbox{  is constant on } \fG \mbox{-orbits}\}. \] Unlike the  case of a group action on $X$, the  groupoid $\fG $ does not in general act on functions.  
\bl \label{gcs}The comorphism $\pi: X \lra \Spec \cO(X)^\fG =Y$ is constant on $\fG$-orbits.\el
\bpf  
This follows from the definition of $\cO(X)^\fG$.
\epf \noi 
For continuous Weyl groupoids the map $\pi$ has some  pleasant properties, see Corollary \ref{sgs6}.

\subsection{Continuous Weyl groupoids and their actions.} \label{cwg}
In \cite{SV2} Sergeev and Veselov associated a certain groupoid $\mathfrak{W}= \mathfrak{W}(\Gd)$,  which they the call Weyl groupoid,  to the  root system
$\Gd $ of a KM  algebra. They also defined a functor $\mathfrak W \lra \mathfrak A(\fh^*)$. We introduce two closely related  groupoids $\mathfrak{W}^c= \mathfrak{W}^c(\Gd)$ and $\mathfrak{W}^c_*= \mathfrak{W}^c_*(\Gd)$. 
\\ \\
We need a preliminary construction, namely
the semi-direct product groupoid $\Gc \ltimes \mathfrak G$.   Let $\mathfrak G$ be a groupoid and  $\Gamma$ a group  acting on $\mathfrak G$  by automorphisms of the corresponding category. In particular, $\Gamma$ acts on the base $\mathfrak B$ of $\mathfrak G$.
Then the  {\it semi-direct product groupoid}
$\Gamma \ltimes \mathfrak G$ has the same base $\mathfrak B$, and 
is defined in \cite{SV2}.
To define the continuous Weyl groupoids, we first introduce some auxilliary groupoids  $\cG^c$  and $\cG^c_*$.

\subsubsection{An action on $\fh^*$.} \label{cwg1} 
In the KM case consider the following groupoid $\mathfrak T^c_{iso}$ with base 
$\Gd_{iso}^c=\Gd_{iso}\ti \ttk.$ 
The non-identity  morphisms  are  $\tau_{\alpha, t}:(\alpha, t) \rightarrow (-\ga, t)$. The group $W$ acts on $\mathfrak T_{iso}^c$ in a natural way: $\alpha \rightarrow w(\alpha),\,
\tau_{\alpha,t} \rightarrow \tau_{w(\alpha),t}$.  Let $\cG^c$  be the semi-direct product groupoid $\cG^c =W \ltimes \mathfrak T_{iso}^c$ with base 
$\Gd_{iso}^c.$  The {\it continuous Weyl groupoid} $\cW^c$ is defined
by 
\be \label{fgh}  \mathfrak{W}^c = W \coprod  \cG^c,\ee     the disjoint union of the group $W$ considered as a groupoid with a single point base $[W]$ and $\cG^c$. 
There is an action of $\mathfrak{W}^c$ on  $\fh^*$, given by a functor from $\mathfrak{W}^c$ to $\mathfrak A(\fh^*)$. Here we adapt \cite{SV2}, Section 9. 
The base point $[W]$ maps to the whole space $\fh^*$.
The base element  $(\alpha,t)$ maps to the hyperplane $\Pi_{\alpha}=\Pi_{-\alpha}$.
The morphism $\tau_{\alpha,t}$ acts via 
\be \label{qyz} \tau_{\alpha,t}(x) = x + t\alpha, \, x \in \Pi_{\alpha}.\ee 
If $ \fg = \fp( n) $ or   $  \fq(n) $, define the groupoid $ \mathfrak T^c$
with base 
$\fB^c=\{(\pm\overline{\ga},t)\mid \ga \in \Gd_0^+, t\in \ttk\} $ and  non-identity morphisms 
$\gt_{\ga,t}:(\overline{\ga} ,t) \lra (-\overline{\ga} ,t)$ and 
$\gt_{-\ga,t}:(-\overline{\ga} ,t) \lra (\overline{\ga} ,t)$ for $\ga \in \Gd_0^+$. Next form the semi-direct product $\cG^c =W \ltimes \mathfrak T^c$ with base $\fB^c$, and define $ \mathfrak{W}^c $ as in \eqref{fgh}.  
There is an action of 
$\mathfrak{W}^c$ on  $ \fh^*$  defined as before, except that in  \eqref{qyz}
  for $\ga \in \Gd_0^+$, 
$\Pi_{\alpha}$ is given   by 
 \eqref{rtak}. The orbits of   
$\mathfrak{W}^c$ on  $ \fh^*$   are described in all cases in Theorem \ref{sgs}.  However, we make no use of this action when $ \fg = \fp( n) $, see Subsection  \ref{pc}. 
\subsubsection{An action on $\mathbb{T}$.} \label{cwg2}
\noi
We need a variant   $\mathfrak{W}_*^c$  of $\mathfrak{W}^c$. In the KM case define $\cG_*^c$  in exactly the same as  $\cG^c$, except that the base of $\cG_*^c$  is 
 $\Gd_{iso}\ti \ttk^*.$ Then set 
$\mathfrak{W}_*^c = W \coprod  \cG_*^c$.  There is an  action of $\mathfrak{W}_*^c $ on
$\mathbb{T}$, 
corresponding to a functor from $\mathfrak{W}_*^c$ to $\mathfrak A(\mathbb{T})$. As before, the base point $[W]$ maps to the whole space $\mathbb{T}$.   The 
 base element  $(\alpha,t)$ maps to the subtorus  $\ttT_\ga$ defined in \eqref{rtaj}.   
The morphism $\tau_{\alpha,t}:\ttT_\ga  \lra \ttT_\ga$ acts via 
\be \label{W_} \tau_{\alpha,t}(x) = c_\alpha(t)x, \, x \in \ttT_\ga.\ee 
 \noi  
If 
$G=P(n) $  or $Q(n)$, The groupoids 
$\cG_*^c $ and  $ \mathfrak{W}_*^c $ are defined by analogy with the groupoids for 
$ \fp( n) $ or   $  \fq(n) $ 
 except that the base of 
 $\cG_*^c$ is $\fB_*^c=
\{(\pm\overline{\ga},t)\mid \ga \in \Gd^+, t\in \ttk^\ti \} .$
There is a functor  $\cW_*^c\lra \fA (\cT)$, such that $W=N_{G_0}(\cT)/\cT$ acts in the obvious  way.   The base elements $(\pm\overline{\ga},t)$ map to $\ttT_\ga = \Ker \tte^{\overline{\ga}}$,  and for  $\ga \in \Gd_0^+$
$\gt_{\ga,t}$ 
acts via \eqref{W_}. 
\bl  \label{qts} If $(\ga,\overline{\gb}) =0$, then $(\gb,\overline{\ga}) =0$.
\el 
\bpf Suppose
 $\ga=\gep_i - \gep_{j},$ and  
$ {\gb}=\gep_k - \gep_{\ell}$ where $i<j, k<\ell$. 
If  $(\ga,\overline{\gb}) =0$, then 
 $0=\gd_{i,k} + \gd_{i,\ell}-\gd_{j,k}-\gd_{j,\ell}.$  If $i=k$, then $j=\ell$ and hence
 $(\gb,\overline{\ga})=\gd_{i,k} - \gd_{i,\ell}+\gd_{j,k}-\gd_{j,\ell}=0.$ Similar arguments for the cases $i=\ell$ and $i\neq k,\ell$ complete the proof.
 \epf 
Equation  \eqref{W_} is a multiplicative  analog of \eqref{qyz}.
If $\ga\in\Gd_{iso}$ is isotropic,  $ \ttT_\ga$ is a codimension one subtorus in $\ttT_\ga$, and by 
\eqref{cwg4}  $\Im \; c_\ga \subseteq \ttT_\ga$. Then $\tau_{\alpha,t}$ translates  $x\in \ttT_\ga$ using the one-parameter subgroup $c_\ga.$ In types P and Q, this holds  since $(\ga,\overline{\ga})=0.$  
\subsubsection{Orbits.} \label{orb}  If $\fG $ is the coutinuous Weyl groupoid as in  Notation \ref{bn}, we describe the orbit   $\fG\gl$ of $\gl\in X$.  
Define the  {\it degree of atypicality} of $\gl$ to be  
$$\atyp \; \gl =  \max\{s|(\gl,A)=0 \mbox{ for some iso-set } A \mbox{ with } |A|  = s \} \quad 
 \mbox{ if  } 
X=\fh^*
,$$ or $$\atyp \; \gl =  \max\{s|\gl\in \bca_{\ga \in A}\ttT_\ga \mbox{ for some iso-set } A \mbox{ with } |A|  = s \} \quad 
 \mbox{ if  } X=\cT.$$
\noi 
We say $\gl $ is {\it typical} if  $\atyp \; \gl =  0.$  Next set
$$E(\gl) = \{\ga\in \Gd_{iso}|\gl\in \Pi_\ga\}, \quad \ttE_\gl=\bcu_{w\in W} w[\bcu_{\ga\in E(\gl)}(\gl + \ttk\ga)]$$ 
 or 
$$E(\gl) = \{\ga\in \Gd_{iso}|\gl\in \ttT_\ga \}, \quad \ttE_\gl=\bcu_{w\in W} w[\bcu_{\ga\in E(\gl)}(\Im c_\ga\gl )],$$ if 
 $X=\fh^*$ or $X=\cT$ respectively.
We need a preliminary result. 
\bl \label{qcs} We have 
$\fG\gl=\ttE_\gl$.\el 
\bpf 
Clearly $\ttE_\gl \subseteq \fG{\gl}.$ We show 
that  for $\mu\in \fG_{(i+1)} \gl,$  
we have $\ttE_\mu =\ttE_\gl $. 
Write 
$\mu  = g\nu$ where $\nu\in \fG_{(i)}\gl$ and  
$g\in \fG_{(\nu)}$. By induction $ \ttE_\nu =\ttE_\gl$. Thus we can assume that $\nu =\gl$. 
Since $\ttE_{\gl}$ is $W$-invariant, we  may  further assume $g=\gt_{\ga,t} $. 
The action of $g$ is given by \eqref{qyz} or 
\eqref{W_}, so since  $g$ is defined at $\gl $, $\gl \in \Pi_\ga$ or $\gl \in \ttT_\ga$.  
Hence  
$\ga\in E(\gl)$ and this   gives the result. 
\epf
\noi 
\br \label{sgs1} {\rm If $\gl\in\fh^*$, the  proof shows that $\cW{\gl} =\bcu_{w\in W} w[\bcu_{\ga\in E(\gl)}(\gl + \Z\ga)],$  which has closure $\cW^c {\gl}$.  
This explains why  our primary interest is in orbits of 
continuous Weyl  groupoids.  
To give a definitive description of these  orbits we need two more Lemmas. }
\er
\bl\label{sgs2}  Let $F(\gl)\subseteq   E(\gl)$  be an iso-set   with $|F(\gl)| =\atyp\;\gl$, and set $G(\gl)
=F(\gl)\cup-F(\gl)$. 
If $\gb\in E(\gl)\bsk F(\gl)$, then for some  reflection $u\in W$ we have 
$u\gl = \gl$ and  $u\gb  = \pm\ga  \in G(\gl).$
\el

\bpf See \cite{Gor1} Equation (14). For convenience we give a short proof. 
In the KM case, we can  argue as follows. There is some $\ga\in F(\gl)$ such that $(\ga,\gb)\neq 0$.  The Jacobi identity applied to the three root vectors $e_{\pm\ga}\in\fg^{\pm\ga}$ and
$ e_\gb\in\fg^{\gb}$ implies that one of $[e_\gb,e_{\pm\ga}]$ is  non-zero.  Replacing ${\ga}$ by ${-\ga}$ if necessary, we can assume ${\gc=\gb}-{\ga}$ is an (even) root.
Then 
from the rank two case, \cite{M} Table 3.4.1  we see that $\gc$-string through $-\gb$ consists of $- \ga$ and $-\gb$. Thus if $u=s_\gc $ is the reflection corresponding to $ \gc$, we have $u \gb=\ga$.  Since $ (\gl,\gc)=0$,  $u$ fixes $\gl .$
For types $P$ and $Q$ we adapt \cite{Gor1}. We have  
 $(\ga,\overline{\gb})\neq 0$ for some $\ga\in F(\gl)$. 
We can assume 
 $\ga = \gep_i - \gep_j, \gb=\pm(\gep_i - \gep_k) \in \Gd_0^+$, where $i<j\neq i, k$. Then 
$$( {\gl},\overline{\gb} ) =   0= ({\gl},\overline{\ga}),$$ 
and the result holds with $\gc = \gep_j - \gep_k$. 
\epf \noi 
Next set 
$$\ttF_\gl=\bcu_{w\in W} w(\gl +\sum_{\ga\in F(\gl)} \ttk\ga) \mbox{ or } \ttF_\gl  = \bcu_{w\in W} w\left(\prod_{\ga\in F(\gl)} \Im c_\ga\right)\gl$$ if 
 $X=\fh^*$ or $X=\cT$ respectively. The definitions do not change if $F(\gl)$ is replaced by $G(\gl)$.
\bl \label{sgs7}   We have $\ttF_\gl\subseteq \fG\gl$.\el\bpf
By \eqref{W_} and   Lemma \ref{qts}, if $\alpha,\gb \in F(\gl)$, $s,t \in \ttk$ and 
$\gl \in \ttT_\ga \cap \ttT_\gb$ then $c_\alpha(t)\gl \in \ttT_\gb$. 
Similar remarks apply to \eqref{qyz} if $\gl \in \Pi_\ga \cap \Pi_\gb$. This implies the statement.\epf
\bt\label{sgs}  We have 
\bi 
\itema  $\fG{\gl} =\ttF_\gl $, 
\itemb $\dim \fG  \gl=  \atyp\;\gl,$
\itemc Every $\fG$-orbit 
is closed.  
\itemd Let $\sim$ be the relation   from Corollary \ref{csn} $($b$)$. 
Then ${\mu\sim\gl}$ iff 
$\fG \mu= \fG\gl$.
\iteme If  $m_\gl =M_\gl \cap \A$ as in Corollary \ref{csn}, then 
$\cV(m_\gl)=\fG  \gl$. 
\ei\et
\bpf 
By  Lemmas \ref{qcs}  and \ref{sgs7}  
$ \ttF_\gl\subseteq  \ttE_\gl =\fG  \gl$. 
We  show   $ \ttE_\gl\subseteq  \ttF_\gl.$  
If $\gb\in E(\gl)\bsk F(\gl)$, let $u $ be as in Lemma \ref{sgs2}.  Then   $u(\gl +\ttk\gb) = 
(\gl + \ttk\ga) $ and $u (c_\gb(t) \gl) = c_\ga( t)\gl$ in the two cases under consideration. This  shows   (a) 
and  (b), (c) are immediate consequences. 
 By definition  $ \sim$ is the smallest equivalence relation satisfying 
the conditions in Corollary \ref{csn} (b).  Hence if ${\mu\sim\gl}$ we have
$\fG\mu= \fG\gl$. Conversely if $\fG\mu= \fG\gl$, then 
 ${\mu\sim\gl} $ by Lemma \ref{qcs},  proving (d).
From Corollary \ref{csn} and  (d) we have 
\be \label{sgs3} m_\gl=m_\mu \mbox{ iff } 
{\mu\sim\gl} \mbox{ iff } \fG \mu= \fG\gl.
\ee
Since $\gl\in\cV(m_\gl)$ this implies $\cV(m_\gl) \supseteq\fG  \gl$. 
For the other inclusion,  if $\mu\in \cV(m_\gl)$, then 
$m_\mu = \cI_\A(\mu ) \supseteq   m_\gl$, and equality must hold. Thus 
$\mu\in\fG  \mu=\fG  \gl,$ proving (e).  \epf
\noi 
We interpret some of our results geometrically.  
The map $\pi$ from Lemma \ref{gcs} restricts to the  map 
 \be \label{sgs5} 
\psi:\Max \cO(X) \lra \Max A. \ee  
as in \eqref{sgs51}.  
 \bc \label{sgs6}  \bi \itema The map in \eqref{sgs5}  is surjective.
\itemb Identifying $X$ with  
$\Max \cO(X)$, the fiber over  $m_\gl$ is   $\fG \gl.$
\ei  
\ec
\bpf (a) is a restatement of the Weak Nullstellensatz, Corollary \ref{csn} (a) and (b) follows from \eqref{sgs3}. 
\epf

\bp \label{qag} Assume Hypothesis \ref{bns} holds. Then 
\bi 
\itema
$\cO(\cT)^{\cW_*^c} =J(G)\ot_\Z \ttk.$
\itemb $S(\fh)^{\cW}=S(\fh)^{\cW^c} = I(\fh) $ if $\fg  \neq \fp(n)$.
\ei
\ep
\bpf  By \eqref{yta}, \eqref{wpn}, \eqref{wqn},   and 
\cite{M22} Lemma 3.1,  
 if $f\in \cO(\cT) $, then $f \in J(G)\ot_\Z \ttk$ iff $f(\lambda) = f(c_\ga({t})\lambda)$ for all relevant roots $\ga$, $t\in\ttk$ and all $\gl \in \ttT_\ga$. Equivalently $f\in \cO(\cT)^{\cW_*^c}.$ This shows (a) and the proof of (b) is similar using \eqref{ltal}. 
The first
 equality in (b) holds because any polynomial that is constant on infinitely many points on a  line is constant on the whole   line.  
\epf   \noi 

\section{The Geometry of Superalgebraic Sets.}\label{csc} 
\subsection{Outline}\label{ocsc} 
Let   $ \fg  $ be a KM Lie superalgebra and let $G$ be the corresponding supergroup as in Lemma \ref{abc}. 
We use Notation \ref{bn}.
We define the S-topology on $\Spec \cO(X) $ by declaring that the S-closed sets are the $\fG$-stable  (Zariski) closed sets. To give an idea of the  geometry of superalgebraic sets, we show that they are  are precisely the the S-closed sets.
We also compute the $S$-closure $V^S$ of an arbitrary closed set $V=V_0=\cV(I)$ with $I$ a $W$-invariant radical ideal in $\cO(X) $. From $V$ we obtain a superalgebraic set which we call  the {\it heart}  and another closed set $V_1=\cV(I_1)$ with $I_1$ a $W$-invariant  radical ideal of $\cO(X) $. The process can be repeated giving new closed sets $V_i$ and new hearts. Then  $V^S$ is the union of all the hearts and also the union  of all the  $V_i$. The hearts are determined by the minimal degrees of atypicality in the Zariski  irreducible components of the $V_i$. A somewhat simplified version of this procedure and some examples are given in \cite{M22}.  
 If the defining equations for $V$ are known, it can be done algorithmically using Groebner bases, see 
 \cite{M22} Example 2.11. The heart is defined in the next subsection, and then $V^S$ is analyzed using projections onto intersections of hyperplanes corresponding to iso-sets.
Let 
$$\atyp \;V = \max\{s| \mbox{ for some } \gl\in V, \atyp \; \gl = s \}.$$
Now suppose $r=\atyp \;V$ and choose an iso-set  $S=\{\gb =\gb_1, \gb_2, \ldots, \gb_r\}$ with  $(\gl,S)=0$ for some $\gl\in V$. Since we can replace $\ga \in S$ by $-\ga$ if necessary, we can assume  $S\subseteq \gO.$ 

\subsection{The heart}\label{hea}
Suppose $I = \bca_{i=1}^\ell P_i$, an irredundant intersection  of prime ideals of $\B$,  where $T\notin P_i$ iff $i\in [k]$. Equivalently $i\in [k]$ iff  $\cV(  P_i)$ contains typical points.  For 
$i\in [k]$, set $Q_i = P_i\cap \A$. Then set
\be \label{ct1} \heart \; V=\heart_0 V = \cV(\bca_{i=1}^k Q_i), \quad V_\reg = \cV(\bca_{i=1}^k P_i). \nn \ee
If $k =0$, set $\heart \; V  = \emptyset.$ 
Note that $\heart \; V$ is superalgebraic, since $ \bca_{i=1}^k Q_i$ is an ideal of $\A$. 
Set $\Vc = \cV(I,T)$.
\bl \label{ct2}  For $i \in [k]$,
\bi\itema $\cV(P_i) \subseteq \cV(Q_i)$ and  
 $V_\reg   \subseteq \heart \; V$. 
\itemb If  $\gl \in \cV(Q_i)$, and $T(\gl)\neq 0$, then $\gl \in \cV(P_i)$. 
\itemc If  $\gl \in \heart \; V$, and $T(\gl)\neq 0$, then $\gl \in V_\reg $.
\itemd $\heart \; V\bsk V_\reg  \subseteq \Vc$.
\iteme If $i>k,$ then  $\cV(P_i)  \subseteq \Vc.$
\ei\el
\bpf(a) is clear since the 
$\cV(P_i)$ and $ \cV(Q_i)$ 
are the irreducible components of  
 $V_\reg $ and $\heart \; V$ respectively.  Next if  $f \in P_i$ we have $T^bf \in  Q_i$ for some $b\ge 0$. So if $\gl \in \cV(Q_i)$, and $T(\gl)\neq 0$ we obtain  $f(\gl)=0,$ proving (b)  and (c) is proved similarly.   By (c), if $\gl\in \heart \; V\bsk V_\reg$, then $T(\gl)=0$ and (d) follows from this. Finally (e) follows since $T\in P_i$ for $i>k$.\epf
\noi Note that 
\be \label{aye}V = V_\reg \cup \bigcup_{i>k} \cV(P_{i}).\ee
\bc \label{5.5} 
\[ V = \heart \; V \cup V^c.\] 
\ec
\bpf The inclusion $\subseteq$  follows from (a), (e) in the Lemma and \eqref{aye}.  The 
opposite inclusion results from  $V^c\subseteq V$, (d)  in the Lemma and \eqref{aye}. 
\epf
 \noi
\subsection{Analysis using projections}\label{anp}
If  $q:X\lra Y$ is a morphism of affine varieties, define the comorphism $q^*:\cO(Y) \lra \cO(X)$ by $(q^*f)(x) =(fq)(x)$ for all $f\in \cO(Y) $ and $x\in X$. It is easy to check the following.
\bl \label{fval} 
Let $i:Y\lra X$ be an embedding of affine varieties and suppose $q:X\lra Y$ is a morphism such that $qi= \id_Y$.  If $I=\ker i^*$, then $$I \cap q^*\cO(Y) =0 \mbox{ and }
I +q^*\cO(Y) =\cO(X). $$ Furthermore, if $V$ is the closed subset of $Y$ defined by the radical ideal $K$ of $\cO(Y)$, then as a closed subset of $X$ we have $V =\cV(q^*K). $
\el

\bexa
\label{em1} {\rm Let $X= \fh^*$ and $X(\gb)=\fh(\gb)^*$.  
Define $q_\gb: \fh^*  \lra \fh(\gb)^*$ 
by restriction $q_\gb(\gl)= \gl |_{\fh(\gb)}$,  and let  $i_\gb: \fh(\gb)^* \lra \fh^*  $ be the map induced by the inclusion $\Gd(\gb)\subset \Gd$.  The map  $q_\gb$ can also be described as follows.  IAN pink notebook
Let $\{\gep_i\}_{i=1}^p$ and $\{h_i\}_{i=1}^p$ be dual bases for $\fh(\gb)^*, \fh(\gb)$.
Extend the $h_i$ to functionals on $\fh^*$ that vanish on the orthogonal complement of 
$ \fh(\gb)^*$ in $\fh^*$, see \eqref{cut}. 
Then for all $\gl\in \fh^*$ we have 
$q_\gb(\gl)= \sum_{i=1}^{p}h_i(\gl) \gep_{i}$
.}
\eexa
\bexa
\label{em2} {\rm 
If $X=\cT$,  let $X(\gb)=\cT(\gb)$ be the subtorus of $X$ generated by $\Im c_\ga$ for $\ga\in \Gd(\gb)$.  
Let $\{\chi_i\}_{i=1}^p$ and $\{c_i\}_{i=1}^p$ be bases for $X(\cT(\gb)), Y(\cT(\gb))$ respectively such that ${\chi_i}c_j(t)= t^{\gd_{i,j}}$.  Extend the $c_i$ to $Y(\cT)$ by analogy with the previous example.  Then define $q_\gb: \cT \lra \cT(\gb) $ by $q_\gb(t)= \prod_{i=1}^p   c_i(\chi_i(t))$.
}\eexa \noi 
Now for $X= \fh^*$ or $\cT$, 
and $\gb$ an isotropic root, define $X(\gb)$ as in the Examples. Define  $p_\gb=i_\gb \ci q_\gb$.  We have 
\be \label{ref} 
 \Pi_\gb =\fh(\gb)^* + \ttk \gb \mbox{ and } \ttT_\gb =\cT(\gb)\ti \Im c_\gb.\ee
\noi 
For  the torus case, $\cT(\gb)\subset \ttT_\gb$ 
by Lemma \ref{cwg3} 
and equality holds because 
$\ttT_\gb$ is  the subtorus generated by $\Im c_\ga$ for $\ga\in \Gd$ with $(\ga,\gb)=0$. 

\bl \label{lm1} \bi\itema  If $\gl\in \Pi_\gb$, then $\gl-p_\gb(\gl)\in \ttk \gb.$
\itemb  If $t\in \ttT_\gb$, then $tp_\gb(t)^{-1}\in \Im c_\gb.$ 
\ei
\el
\bpf We prove (a).  The proof of (b) is similar.  By \eqref{ref} 
we can write $\gl=  \gl_0 + c\gb$ with 
$\gl_0 \in \fh(\gb)^*$, $c\in \ttk.$ Since 
$i_\gb\ci q_\gb$ is the identity on $\fh(\gb)^*$ and $q_\gb(\gb)=0$, we obtain the result.
\epf \noi 
The action of $W$ on $\cO(X)$ is given by $(wf)(x) = f(w^{-1}x)$ for $w\in W, f\in \cO(X)$ and $x  \in X$.  Then $w  \Gd(\gb)= \Gd(w\gb)$, and we have 
$w X(\gb)= X(w\gb).$ Since $W$ acts transitively
on $\gO$ we may assume $\gs(w\gb)=w\gs(\gb)$ in \eqref{cut}.  Then
\be \label{cot} 
w\ci q_\gb(\gl)=q_{w\gb}(w\gl).\ee
We have an embedding $q_\gb^*:\cO(X(\gb))\lra \cO(X)$
as in Lemma \ref{fval} with $Y=X(\gb)$ and $q=q_\gb.$ In what follows we suppress
$q_\gb^*.$  
\\ \\
Set $X_\gb= \Pi_\gb$ or $\Pi_\gb=\ttT_\gb$ 
in  case (i) or (ii) of \ref{bn} holds
respectively
For $\gb\in \Gd_{iso}$, set  $V_\gb =V\cap X_\gb$. 
Since $T(\gl)=0$ for $\gl \in\Vc $, we have
\be \label{ct3}\Vc = \bigcup_{\gb\in \gO}
 V_{\gb }.\nn\ee \noi 
Next we study  the sets $ V_{\gb}$ in more detail.  
If $\gb = \gb_1$, set $ U_{\gb}= q_{\gb}(V_{\gb})\subseteq X({\gb})$ and 
\by \label{e41}
K_{\gb}&=&  \cI(V_{\gb}) \cap \cO(X({\gb})), \nn\ey  
We denote the closure of $U\subseteq X$ by $\overline{U}$. 
\noi  
\bl \label{zyo} For all ${\gb\in\gO } $, in case $($i$)$ $V^S$ contains the sets 
\be \label{xyp} Y_{\gb}:= \overline{U_{\gb} +\ttk\gb} =\overline{U_{\gb}} +\ttk\gb.
\ee 
\el 
\bpf If $\gl\in V_{{\gb }}$, then by Lemma \ref{lm1}, $\gl-p_{\gb}(\gl) =a{\gb }$ for some $a\in \ttk$. 
By $\cW^c$-invariance, $V^S$ contains $ \gl +(t-a){\gb } = p_{\gb}(\gl) + t{\gb }$ for all $t\in \ttk.$ For the second equality in \eqref{xyp}, note that $\overline{U_{\gb}} \cap \ttk\gb=0.$
\epf \noi 
Similarly in case (ii) $V^S$ contains the sets 
$$ Y_{\gb}:= \overline{U_{\gb} \ti \Im c_\gb} =\overline{U_{\gb}} \ti \Im c_\gb.$$ 
Next if $X= \fh^*$ set $z_{\gb} = h_{\gs(\gb)}$, and   if $X=\cT$
set $z_{\gb}= \tte^{\gs(\gb)}-1$.  Then 
let $ L_\gb$ be the ideal of $\cO(X)$ generated by $K_\gb$ and $z_{\gb}$.  
\bl \label{CL2}
We have $\overline{U}_{\gb}=\cV(K_\gb)$ and $Y_{\gb}  =\cV(L_{\gb})$.
\el \noi 
\bpf  The first statement follows  by \cite{CLO}, Theorem 3.2.3, page 122 and the second is an immediate consequence.
\epf  \noi 
The ideal $K_{\gb}$ is called an {\it elimination ideal}.
Next set  $V_1=\bcu_{\gb\in\gO} Y_\gb$. Then  $V_1=\cV (I_1)$ where $I_1=\bca_{\gb\in\gO} L_\gb$.  Because  $I=\cI(V)$  is a radical ideal, so too are $K_{\gb}, L_{\gb}$  and $I_1$.  By Lemma \ref{zyo} we have
\bp \label{DL} 
The set $V^S$ contains $\heart \; V\cup   V_1$.  So, $V^S=\heart \; V\cup   V_1^S=  V\cup   V_1^S.$\ep \noi 
The next result shows that $V_1$ and $I_1$ are $W$-invariant. 
\bl \label{AL}  If $V$ is $W$-invariant, then  for $w\in W,$ we have 
$w(V_\gb)=V_{w\gb}$
\bi \itema $w(V_\gb)=V_{w\gb}$, 
$w(U_{\gb})= U_{w\gb} $ and $w(Y_\gb)=Y_{w\gb}$
\itemb $w(K_\gb\cO(X))=K_{w\gb}\cO(X)$ and  $w(L_\gb)=L_{w\gb}.$ \ei
\el 
\bpf The first statement in (a) follows from $W$-invariance of the bilinear form $(\;,\;)$, the second from \eqref{cot} and the third is an immediate consequence.  
By Lemma \ref{fval}, the closed sets defined by  $w(K_\gb\cO(X))$ and $K_{w\gb}\cO(X)$ are $w(\overline{U}_{\gb})$ and $ \overline{U}_{w\gb} $.  Since these are equal, the defining ideals are equal.  This gives the first part in (b) and the second also follows.
\epf \noi
Now we repeat this process. 
Set  $A=A(q)=\{\gb_1,\ldots,\gb_q \} $ and 
$V_A= V\cap \bca_{\gb\in A} \Pi_\gb$ in case (i) or 
$V_A= V\cap \bca_{\gb\in A} \ttT_\gb$  in case (ii). 
Define  $p_{A}: X\lra X$ by $p_{A} = p_{\gb_q}\ci \ldots \ci p_{\gb_1}.$ 
\noi Next, if  $X_A=\bca_{\gb\in A} X(\gb)$, then
$\cO(X_A)=\bca_{\gb\in A} \cO(X(\gb)).$
Consider the elimination ideal 
\by \label{e5}
K_{A}&=&  \cI(V_A) \cap \cO(X_A).\nn \ey
Let $ L_A$ be the ideal of $\cO(X )$ generated by $K_A$ and $z_{\gb_1},\ldots, z_{\gb_q}.$ 
If $ U_{A}= p_{\ga}(V_{A})$, then  $\overline{U}_{A}=\cV(K_A)$. 
\bl \label{lm5} For all $ q\le r$, $V^S$ contains  the sets 
\be \label{xyq} \{ p_{A(q)}(\gl) + \sum_{\gb\in A(q)}\ttk  {\gb }|\gl\in V_{A(q)} \}.\nn
\ee
\be \mbox{  or } \{ p_{A(q)}(\gl) 
\prod_{\gb\in A(q)} \Im c_\ga
|\gl\in V_{A(q)} \} \mbox{  in case } (ii).\nn\ee
Hence  $V^S$ also contains  the closures $  Y_{A(q)}=\cV(L_{A(q)})$ of these sets. 
\el
\bpf This follows as in the proof of Lemmas \ref{zyo} and \ref{CL2}, by $\fG$-invariance.\epf \noi 
We define $V_q = \bcu_{w\in W} wY_{A(q)}$ and
$I_q = \bca_{w\in W} wL_{A(q)}$. Then 
$V_q = \cV(I_q)$. 
By construction if $\gl \in V_q$ then  $\atyp \; \gl \ge q.$
Suppose $I_q = \bca_{i=1}^\ell P_i$, an irredundant intersection  of prime ideals,  and suppose  $i\in [k]$ iff  $\cV(  P_i)$ contain a point $\gl$ with $\atyp \; \gl = q.$  For 
$i\in [k]$, set $Q_i = P_i\cap A$. Then set
$ \heart_q V = \cV(\bca_{i=1}^k Q_i).$ 
If $k =0$, set $\heart_q V  = \emptyset.$ 

\bt \label{xya} 
The S-closure of $V$ is 
$V^S = \bcu_{i=0}^r\heart_i V= \bcu_{i=0}^r V_i.$  The set $V^S$ is superalgebraic. 
\et
\bpf As in Proposition \ref{DL} if $q<r$, we have $V_{q}^S=\heart_q  V\cup   V_{q+1}^S =V_q  \cup   V_{q+1}^S$ and $V_{r}^S =\heart_r  V$. This gives the equalities. Each $\heart_i V$ is superalgebraic, so the result follows.
\epf
\bc\label{ass} The superalgebraic sets are exactly the closed sets that are invariant under the Weyl groupoid  $\fG$.\ec 

\bpf If  $V$ is closed and 
$\fG$-invariant, then  
$V= V^S$ is superalgebraic.  Conversely a superalgebraic set has the form $V = \cV(I)$, with $I$ an ideal of $\A=\cO(X)^\fG$. Since elements of $\A $ are constant on orbits, $V$ is a union of $\fG$-orbits.
\epf

\end{document}